\documentclass[11pt,a4paper]{article}
\usepackage[a4paper,margin=1in]{geometry}
\usepackage[T1]{fontenc}
\usepackage{lmodern}
\usepackage{amsmath,amssymb,amsthm,mathtools}
\usepackage{microtype,booktabs,indentfirst}
\usepackage{mathrsfs}
\usepackage{xcolor}
\usepackage[colorlinks=true,linkcolor=blue,citecolor=blue,urlcolor=blue]{hyperref}

\allowdisplaybreaks[2]

\numberwithin{equation}{section}
\newtheorem{theorem}{Theorem}[section]
\newtheorem{proposition}[theorem]{Proposition}
\newtheorem{lemma}[theorem]{Lemma}
\newtheorem{corollary}[theorem]{Corollary}

\theoremstyle{definition}

\newtheorem{remark}[theorem]{Remark}
\newtheorem*{problemfour}{Open Problem 11.4}
\newtheorem*{problemfive}{Open Problem 11.5}

\title{Local minimizers of higher Dirichlet eigenvalues\thanks{{Research supported by the National Natural Science Foundation of China} (No. 12371110).}}

\author{
	Guowei Dai\thanks{School of Mathematical Sciences,
		Dalian University of Technology, Dalian 116024, P.R. China.
		Email: daiguowei@dlut.edu.cn.}
	\and
	Lan Zhao\thanks{Corresponding author. School of Mathematical Sciences,
		Dalian University of Technology, Dalian 116024, P.R. China.
		Email: zhaolan@mail.dlut.edu.cn.}}

\date{}
\begin{document}
\maketitle
\begin{abstract}
For every dimension $d\ge2$, we prove that the ball is a local minimizer of the $(d+1)$-th Dirichlet eigenvalue under the volume constraint. Moreover, we establish a quantitative estimate. As a corollary, we also establish quantitative local minimality of the associated unequal two-ball structure for the $(d+2)$-th eigenvalue under smooth perturbations preserving the two component structure and total volume. Our approach combines a matrix representation of the perturbed
eigenvalues, spherical harmonic analysis, and uniform estimates
for the higher-order remainder.
\end{abstract}


\section{Introduction}

\paragraph{Background and motivation.}
Let $\Omega\subset\mathbb{R}^d$, $d\ge2$, be an open set of finite
positive volume $\left|\Omega\right|$. We consider the Dirichlet
Laplacian $-\Delta$ on $\Omega$. It is well known that its spectrum
is discrete and consists of eigenvalues
\[
0<\lambda_1\left(\Omega\right)
\le\lambda_2\left(\Omega\right)
\le\cdots
\le\lambda_n\left(\Omega\right)
\le\cdots.
\]
listed in nondecreasing order and counted with multiplicity. The study of spectral shape optimization is a classical topic in spectral geometry and goes back to Rayleigh's question about the lowest vibration frequency of a membrane of prescribed area~\cite{Rayleigh}. For a membrane fixed along its boundary, Rayleigh proposed that the circular membrane has the lowest fundamental frequency. The problem is simple to state, but the
conjecture was proved almost half a century later. Faber~\cite{Faber} proved the planar result in 1923, and Krahn~\cite{Krahn25,Krahn26} established the corresponding inequality in arbitrary dimension. Thus, the Faber--Krahn theorem tells us that the ball is the minimizer of the first Dirichlet eigenvalue among domains of equal volume. However, the second Dirichlet eigenvalue behaves differently. By the Krahn--Szeg\H{o} inequality, the minimizer of $\lambda_2$ consists of two disjoint balls of equal volume \cite{Henrot,Krahn26,PS}. Ashbaugh and Benguria \cite{AB} subsequently showed that the ball maximizes the ratio $\lambda_2/\lambda_1$. Beyond these sharp extremal configurations, quantitative bounds relate the eigenvalue gap to the domain’s deviation from the optimal geometry. Notably, Brasco, De Philippis, and Velichkov \cite{BDPV} established the sharp quantitative version of the Faber--Krahn inequality.

Changing the boundary conditions gives rise to new phenomena and can also render the problem more complex. For the Neumann problem, Szeg\H{o}~\cite{Szego} proved in 1954 that the disk maximizes the first positive Neumann eigenvalue among bounded
simply connected planar Lipschitz domains of fixed area. Two years later, Weinberger~\cite{Weinberger}
extended this result to all dimensions $d\ge 2$ without assuming simple
connectedness. These results show that both the spectral index and the boundary condition affect the optimal shape and the direction of the inequality. The case of Robin boundary conditions is more delicate, and shape optimization problem for the first eigenvalue remains unresolved in full generality. The Robin boundary condition takes the form \(\partial_{\nu}u+\beta u=0\) on a smooth boundary, with \(\nu\) the outward unit normal. For each fixed $\beta>0$, Bossel~\cite{Bossel} proved the planar Robin Faber--Krahn inequality in 1986, and Daners~\cite{DanersRobin} proved it in every dimension in 2006. The ball therefore remains a minimizer of the first eigenvalue in this setting. Kennedy~\cite{KennedyRobin} subsequently established the corresponding result for the second eigenvalue. As expected, two equal balls minimize the second Robin eigenvalue under a volume constraint. However, a striking and unexpected phenomenon arises for negative values of the parameter. Freitas and Krej\v{c}i\v{r}\'{i}k~\cite{FKRobin} disproved the conjecture that the ball always maximizes the first Robin
eigenvalue. More precisely, they proved that, for sufficiently
large negative values of the Robin parameter, a spherical shell
has a larger first Robin eigenvalue than a ball of
the same volume. In dimension two, they also proved that the disk maximizes the first Robin eigenvalue among bounded \(C^2\) domains of fixed area whenever \(\beta<0\) and \(|\beta|\) is sufficiently small. The Robin problem also shows why a simple geometric intuition cannot be assumed for higher eigenvalues. Kennedy~\cite{KennedyHigher} proved that, in several higher eigenvalue cases, a minimizing shape cannot be independent of the positive boundary parameter. Thus, the optimality of balls for the first eigenvalues cannot be extended directly to higher spectral optimization problems. These examples also show that the shape optimization problem for eigenvalues is more subtle than one might expect. 

We next return to the Dirichlet problem and address the question of when a ball remains optimal for higher eigenvalues. Several important results in this direction have been established. Bucur~\cite{Bucur} proved that for every $k\geq 1$, the $k$th Dirichlet eigenvalue
attains its minimum among quasi-open sets under a volume constraint and that every
minimizer is bounded. Numerical results by Oudet~\cite{Oudet} and Antunes and Freitas~\cite{AF} provide candidate minimizers. These works give a framework for minimization, but they do not identify the optimal shape for the higher eigenvalues considered below. Antunes and Oudet formulated two such problems in Chapter 11 of \cite{AO}. We here present them in a common framework. Let $B$ be the unit ball in $\mathbb{R}^d$ and $\omega_d=\left|B\right|$. We write $B_R$ for the ball of radius $R$ centered at the origin. The notation $j_{\mu,k}$ denotes the $k$-th positive zero of the Bessel function $J_\mu$~\cite{DLMF,Watson}.

\begin{problemfour}[See~\cite{AO}]
The $d$-dimensional ball minimizes $\lambda_{d+1}$ among all $d$-dimensional domains of volume $\mathcal V$.
\end{problemfour}

\begin{problemfive}[See~\cite{AO}]
Among all $d$-dimensional domains of volume $\mathcal V$, a minimizer of $\lambda_{d+2}$ is the disjoint union of two balls whose larger and smaller radii satisfy
\[
\frac{R_1}{R_2}=\frac{j_{d/2,1}}{j_{d/2-1,1}}.
\]
\end{problemfive}

In dimension two, these problems concern the third and fourth Dirichlet eigenvalues. However, despite their apparent simplicity, their global statements remain open \cite{AO,Henrot}. The corresponding local questions have been studied in low dimensions.
Wolf and Keller~\cite{WK} proved local minimality of the disk for
$\lambda_3$. Using the eigenvalue expansions of Wolf and Keller~\cite{WK},
Osting and Kao~\cite{OK} studied convex combinations of
the first three Dirichlet eigenvalues near the disk. When the first-order
splitting vanishes, their formula gives the second-order term in the
expansion of $(\lambda_2+\lambda_3)/2$ under area-preserving deformations. 

\paragraph{Main results.} In this paper, we study the local case of these questions. Our aim is to prove a positive second-order estimate at the ball in every dimension \(d\ge 2\) and to establish a quantitative inequality. Let $\mathbb{S}^{d-1}=\partial B$ and $s_d=|\mathbb{S}^{d-1}|=d\omega_d$. For a real-valued function
$h\in C^{2,\alpha}(\mathbb{S}^{d-1})$, we denote by $h_n$
its spherical harmonic component of degree $n$. We use $\left\lVert\cdot\right\rVert_2$ for the norm in $L^2(\mathbb{S}^{d-1})$ and the following equivalent norm on the boundary Sobolev space $H^{1/2}(\mathbb{S}^{d-1})$:
\begin{equation}\label{eq:hnorm}
\left\lVert h\right\rVert_{1/2,*}^2
=\sum_{n=0}^{\infty}\left(1+n\right)\left\lVert h_n\right\rVert_2^2.
\end{equation}
By orthogonality of the spherical harmonics and the bound $1+n\ge1$,
we have $\|h\|_2\le\|h\|_{1/2,*}$. By the continuous embedding
$C^{2,\alpha}(\mathbb{S}^{d-1})
\hookrightarrow H^{1/2}(\mathbb{S}^{d-1})$, we have
$\left\lVert h\right\rVert_{1/2,*}
\le C\left\lVert h\right\rVert_{C^{2,\alpha}}$,
where $C>0$ depends only on $d$ and $\alpha$. For a sufficiently small function $h$ of zero mean, we construct an exactly volume-preserving radial deformation by
\begin{equation}\label{eq:chart}
\Omega_h
=\left\{r\theta:\theta\in\mathbb{S}^{d-1},\ 0\le r<\zeta\left(h\right)\left(1+h\left(\theta\right)\right)\right\},
\end{equation}
where
\begin{equation}\label{zeta(h)}
\zeta\left(h\right)
=\left(\frac{1}{s_d}\int_{\mathbb{S}^{d-1}}\left(1+h\right)^d\,\mathrm{d} S\right)^{-1/d}.	
\end{equation}
Thus $\left|\Omega_h\right|=\omega_d$. We denote the common eigenvalue of the first angular cluster by
\begin{equation}
\lambda_*:=\lambda_2\left(B\right)=\cdots=\lambda_{d+1}\left(B\right)=j_{d/2,1}^2.
\end{equation}
These eigenvalues form an isolated cluster of multiplicity $d$. Here, an \emph{isolated cluster} is a group of eigenvalues, counted with
multiplicity, separated from the rest of the spectrum by positive
gaps. Under sufficiently small perturbations, these eigenvalues
may split, but they remain an isolated cluster of total multiplicity $d$. Our first theorem is as follows. The scalar coefficients $A_{d,n}$ depend only on the dimension and the harmonic degree; their explicit formula is given in~\eqref{eq:coefficients} in Section~\ref{sec:prelim}.

\begin{theorem}\label{thm:kernel}
Let $h\in C^\infty(\mathbb{S}^{d-1})$ be real and satisfy $h_0=h_1=h_2=0$. Then
\begin{equation}\label{eq:maintrace}
\frac1d\sum_{j=2}^{d+1}\lambda_j\left(\Omega_{\varepsilon h}\right)
=\lambda_*+\frac{2\lambda_*}{s_d}\varepsilon^2
\sum_{n\ge3}A_{d,n}\left\lVert h_n\right\rVert_2^2+o_h\left(\varepsilon^2\right).
\end{equation}
For every fixed $d\ge2$, there is a constant $c_d>0$ such that
\begin{equation}\label{eq:coercivity}
A_{d,n}\ge c_d\left(1+n\right)\qquad n\ge3.
\end{equation}
Consequently,
\[
\liminf_{\varepsilon\to0}
\frac{\lambda_{d+1}\left(\Omega_{\varepsilon h}\right)-\lambda_*}{\varepsilon^2}
\ge\frac{2\lambda_*c_d}{s_d}\left\lVert h\right\rVert_{1/2,*}^2.
\]
Here $o_h(\varepsilon^2)/\varepsilon^2\to0$ for each fixed $h$.
\end{theorem}

Theorem~\ref{thm:kernel} implies that, for each fixed nonzero $h$
satisfying $h_0=h_1=h_2=0$, there exists $\varepsilon_h>0$ such that
\[
\lambda_{d+1}\left(\Omega_{\varepsilon h}\right)>\lambda_*
\qquad\text{for }0<|\varepsilon|<\varepsilon_h.
\]
However, $\varepsilon_h$ may depend on the deformation direction $h$,
so this result does not establish local minimality in a common
neighborhood of the ball. To overcome this difficulty, we establish
a uniform estimate that is valid for all sufficiently small
deformations and also incorporates the degree-two component. We set
\[
\mathscr D\left(h\right)
=\left\lVert h_2\right\rVert_{1/2,*}+\left\lVert h_{\ge3}\right\rVert_{1/2,*}^2,
\qquad h_{\ge3}=\sum_{n\ge3}h_n.
\]

\begin{theorem}\label{thm:uniform}
For each $d\ge2$ and $0<\alpha<1$, there are constants $\eta,c,C>0$, depending only on $d$ and $\alpha$, such that every real $h\in C^{2,\alpha}(\mathbb{S}^{d-1})$ satisfying $h_0=h_1=0$ and $\left\lVert h\right\rVert_{C^{2,\alpha}}<\eta$ obeys
\begin{equation}\label{eq:mainuniform}
c\mathscr D\left(h\right)
\le\lambda_{d+1}\left(\Omega_h\right)-\lambda_*
\le C\mathscr D\left(h\right).
\end{equation}
After decreasing $\eta$ and changing $c$ if necessary, we also have
\begin{equation}\label{eq:quadraticlocal}
\lambda_{d+1}\left(\Omega_h\right)-\lambda_*
\ge c\left\lVert h\right\rVert_{1/2,*}^2.
\end{equation}
\end{theorem}

\paragraph{Difficulties, ideas and contributions.}
The difficulty we have to address is to obtain estimates that are
uniform over all harmonic degrees. As mentioned earlier, positivity along each fixed
direction alone does not imply local minimality in a common
neighborhood. To overcome this difficulty, we establish a uniform estimate on the remainder. Using the formula for change of variables and Kato's spectral projection construction~\cite{HP,Kato,SZ}, together with
the two-norm approach of Dambrine and Lamboley~\cite{DL}, we bound the remainder by $\| h\|_{C^{2,\alpha}}\| h\|_{1/2,*}^2$. This allows us to absorb the remainder
uniformly and combine the first- and second-order estimates to
prove Theorem~\ref{thm:uniform}. 

A key ingredient is Lemma~\ref{lem:cluster}, in which we establish
a matrix decomposition and uniform estimates for general deformations
of an arbitrary smooth domain $\Omega$. This lemma plays an important
role in the proof of Theorem~\ref{thm:uniform}. Moreover, since the eigenvalue is not simple, the classical Hadamard
formula for simple eigenvalues is not applicable. Therefore, we 
use the corresponding result for multiple eigenvalues. Indeed, the perturbed eigenvalues in the cluster can be represented as the eigenvalues
of a real symmetric matrix. This enables us to use matrix analysis and
linear operator techniques on finite-dimensional spaces in the proof
of Lemma~\ref{lem:cluster}.

Thus Theorem~\ref{thm:uniform} gives an affirmative answer to the local
case of Open Problem~11.4 under small $C^{2,\alpha}$
perturbations when volume is fixed. Section~\ref{sec:uniform} subsequently presents two
corollaries of this result. Corollary~\ref{cor:asymmetry} gives a lower bound for the eigenvalue difference in terms of the squared Fraenkel asymmetry, which quantifies the deviation from the ball. Corollary~\ref{cor:twoball}
establishes quantitative local minimality of the unequal
two-ball configuration for $\lambda_{d+2}$, controlling both
changes in the component shapes and changes in their volume
split. Thus it provides an affirmative answer to the local
version of Open Problem~11.5 under small $C^{2,\alpha}$
perturbations preserving the two components and total volume.

\paragraph{The organization of the paper.} The paper is organized as follows.
Section~\ref{sec:prelim} describes the radial deformation, the corresponding perturbation of the spectral subspace, and the boundary operator used in the proof.
Section~\ref{sec:trace} proves Theorem~\ref{thm:kernel}.
Section~\ref{sec:uniform} proves Theorem~\ref{thm:uniform}
and Corollaries~\ref{cor:asymmetry} and~\ref{cor:twoball}.
Appendix~\ref{app:bessel} contains the Bessel zero and boundary operator estimates. Appendix~\ref{app:coefficientweights} bounds the Taylor coefficients used in Lemma~\ref{lem:cluster}.

\section{Preliminaries}\label{sec:prelim}
We first describe the volume-preserving radial deformations of the ball.

\subsection{Volume-preserving deformations}
We denote by $\mathcal H_n$ the space of restrictions to $\mathbb{S}^{d-1}$ of homogeneous harmonic polynomials of degree $n$. The orthogonal projection in $L^2(\mathbb{S}^{d-1})$ onto this space is denoted by $\Pi_n$, so $h_n=\Pi_n h$. The Laplace--Beltrami operator $\Delta_{\mathbb{S}^{d-1}}$ is the tangential Laplacian on the sphere. Spherical harmonics satisfy the eigenvalue equation
for the Laplace--Beltrami operator on the sphere~\cite{ABR}:
\[
-\Delta_{\mathbb{S}^{d-1}}Y
=n\left(n+d-2\right)Y
\qquad Y\in\mathcal{H}_n.
\]
Moreover, they give the orthogonal decomposition $L^2(\mathbb{S}^{d-1})
=\widehat{\bigoplus}_{n=0}^{\infty}\mathcal{H}_n$, where the hat denotes completion in $L^2$.

The norm~\eqref{eq:hnorm} is equivalent to the usual Sobolev trace norm~\cite{McLean}. Finite sums of spherical harmonics are dense in this norm. We denote the dual of $H^{1/2}(\mathbb{S}^{d-1})$ by \(H^{-1/2}(\mathbb{S}^{d-1})\). The pairing \(\left\langle f,g\right\rangle\) coincides with \(\int_{\mathbb{S}^{d-1}}fg\,\mathrm{d} S\) when both functions are smooth and real. The next elementary expansion concerns the particular chart~\eqref{eq:chart}.

\begin{lemma}\label{lem:volume}
	\textnormal{(i)} Let $\Omega_\varepsilon$ be a smooth
	volume-preserving radial deformation of the unit ball with boundary $\partial\Omega_\varepsilon
	=\big\{\rho_\varepsilon(\theta)\theta:
	\theta\in\mathbb S^{d-1}\big\}$, where the boundary radius has the expansion $\rho_\varepsilon
	=1+\varepsilon h+\varepsilon^2g+O(\varepsilon^3)$. Then the coefficients satisfy
	\begin{equation}\label{eq:volumeconstraint}
		\int_{\mathbb S^{d-1}}h\,\mathrm{d}S=0,
		\qquad
		\int_{\mathbb S^{d-1}}g\,\mathrm{d}S
		=-\frac{d-1}{2}\int_{\mathbb S^{d-1}}h^2\,\mathrm{d}S.
	\end{equation}
	\textnormal{(ii)} For every real zero-mean function
	$h\in C^{2,\alpha}(\mathbb S^{d-1})$ with sufficiently small
	$C^{2,\alpha}$ norm, we have the uniform estimate
	\begin{equation}\label{eq:zetaremainder}
		\zeta\left(h\right)
		=1-\frac{d-1}{2s_d}\left\lVert h\right\rVert_2^2
		+r_\zeta\left(h\right),
		\qquad
		\left|r_\zeta\left(h\right)\right|
		\le C\left\lVert h\right\rVert_{L^\infty}
		\left\lVert h\right\rVert_2^2.
	\end{equation}
	In particular, for each fixed smooth real function $h$ of zero mean
	and sufficiently small $|\varepsilon|$, the domain
	$\Omega_{\varepsilon h}$ in~\eqref{eq:chart} has volume $\omega_d$
	and boundary
	$\partial\Omega_{\varepsilon h}
	=\big\{\rho_\varepsilon(\theta)\theta:
	\theta\in\mathbb S^{d-1}\big\}$.
	Its boundary radius has the expansion
	\begin{equation}\label{eq:radialpath}
		\rho_\varepsilon
		=\zeta\left(\varepsilon h\right)\left(1+\varepsilon h\right)
		=1+\varepsilon h+\varepsilon^2g+O_h(\varepsilon^3),
		\qquad
		g=-\frac{d-1}{2s_d}\left\lVert h\right\rVert_2^2.
	\end{equation}
\end{lemma}

\begin{proof}
	\textnormal{(i)} By direct integration, we have
	\begin{align*}
		\left|\Omega_\varepsilon\right|
		&=\frac1d\int_{\mathbb{S}^{d-1}}
		\rho_\varepsilon^d\,\mathrm{d}S\\
		&=\omega_d
		+\varepsilon\int_{\mathbb{S}^{d-1}}h\,\mathrm{d}S
		+\varepsilon^2\int_{\mathbb{S}^{d-1}}
		\left(g+\frac{d-1}{2}h^2\right)\,\mathrm{d}S
		+O(\varepsilon^3).
	\end{align*}
	Since $\left|\Omega_\varepsilon\right|=\omega_d$ for all sufficiently
	small $|\varepsilon|$, the coefficients of $\varepsilon$ and
	$\varepsilon^2$ must vanish. This proves~\eqref{eq:volumeconstraint}.
	
	\textnormal{(ii)} Since $\left\lVert h\right\rVert_{C^{2,\alpha}}$
	is sufficiently small, we may assume that
	$\left\lVert h\right\rVert_{L^\infty}\le1/2$.
	Since $d$ is an integer and
	$\int_{\mathbb S^{d-1}}h\,\mathrm{d}S=0$,
	the binomial formula yields
	\[
	\frac1{s_d}\int_{\mathbb{S}^{d-1}}
	\left(1+h\right)^d\,\mathrm{d}S
	=1+\frac{d\left(d-1\right)}{2s_d}
	\left\lVert h\right\rVert_2^2+r_0,
	\qquad
	\left|r_0\right|
	\le C\left\lVert h\right\rVert_{L^\infty}
	\left\lVert h\right\rVert_2^2.
	\]
	We set
	$q_0=\frac{d\left(d-1\right)}{2s_d}
	\left\lVert h\right\rVert_2^2+r_0$.
	By~\eqref{zeta(h)}, Taylor's formula for the scalar function
	$q\mapsto(1+q)^{-1/d}$ gives
	\[
	\zeta\left(h\right)
	=(1+q_0)^{-1/d}
	=1-\frac{q_0}{d}+r_1,
	\qquad
	\left|r_1\right|\le Cq_0^2.
	\]
	Moreover,
	\[
	q_0^2
	\le C\left\lVert h\right\rVert_2^4
	\le Cs_d\left\lVert h\right\rVert_{L^\infty}^2
	\left\lVert h\right\rVert_2^2
	\le C\left\lVert h\right\rVert_{L^\infty}
	\left\lVert h\right\rVert_2^2.
	\]
	We obtain~\eqref{eq:zetaremainder} after combining the two remainders.
	By the definition of $\zeta$, the domain $\Omega_{\varepsilon h}$
	has volume $\omega_d$. Replacing $h$ by $\varepsilon h$
	in~\eqref{eq:zetaremainder} and multiplying by $1+\varepsilon h$
	proves~\eqref{eq:radialpath}.
\end{proof}

\begin{remark}
Every smooth domain of volume $\omega_d$ sufficiently close to the unit ball can, after a translation, be represented by \eqref{eq:chart} with $h_0=h_1=0$. This translation leaves the Dirichlet eigenvalues and the volume unchanged; see also~\cite{DL,HP}.
\end{remark}

\subsection{Spectral properties and perturbation theory}\label{subsec:perturbation}
We use the standard Sobolev space $H^1(\Omega)$. The space $H_0^1(\Omega)$ is the closure of the smooth compactly supported functions in the $H^1$ norm. The Dirichlet Laplacian is the operator associated with the quadratic form
\[
u\longmapsto\int_\Omega\left|\nabla u\right|^2\,\mathrm{d} x,
\qquad u\in H_0^1\left(\Omega\right).
\]
For a bounded smooth domain $\Omega$, the embedding $H_0^1(\Omega)\hookrightarrow L^2(\Omega)$ is compact. The form therefore defines a positive self-adjoint operator having compact resolvent~\cite{Davies}. Its spectrum
is discrete and consists of eigenvalues. For $r>0$ and $y\in\mathbb{R}^d$, the Dirichlet eigenvalues and
volume are invariant under translation and obey the following
scaling laws:
\begin{equation}\label{scaling-laws}
	\lambda_k\left(r\Omega+y\right)=r^{-2}\lambda_k\left(\Omega\right),
	\qquad
	\left|r\Omega+y\right|=r^d\left|\Omega\right|.
\end{equation}
On a disjoint union of open components, the Dirichlet operator is
the direct sum of the component operators, and its spectrum is
the union of their spectra, arranged in nondecreasing order and
counted with multiplicity~\cite{Davies}.

In polar coordinates $x=r\theta$, the Euclidean Laplacian is
\[
\Delta=\partial_{rr}+\frac{d-1}{r}\partial_r+\frac{1}{r^2}\Delta_{\mathbb{S}^{d-1}}.
\]
The regular separated radial solutions are proportional to $r^{1-d/2}J_{\ell+d/2-1}(\sqrt\lambda\,r), \ell\ge0$. Thus the ball eigenvalues are $j_{\ell+d/2-1,k}^2$. The monotonicity of the first zero and the interlacing of consecutive orders in~\cite[Section~10.21]{DLMF} imply
\[
j_{d/2-1,1}<j_{d/2,1}<j_{d/2-1,2},
\qquad
j_{d/2,1}<j_{d/2+1,1}.
\]
It follows that
\begin{equation}
\lambda_1\left(B\right)=j_{d/2-1,1}^2<\lambda_*<\lambda_{d+2}\left(B\right).
\end{equation}
These gaps will keep the labels $2,\ldots,d+1$ fixed during a sufficiently small deformation.

Throughout this paper, we denote by $\|\cdot\|$ the operator norm
of a bounded linear operator on a Hilbert space. We denote the identity operator on a space $X$ by $\operatorname{Id}_X$
and the $n\times n$ identity matrix by $I_n$.
For a matrix $A=(a_{ik})$, we write $\det A$ for its determinant. In particular,
for a matrix $A\in\mathbb R^{m\times m}$, we have
\[
\|A\|=\sup_{|x|=1}|Ax|.
\]
We denote by $\|\cdot\|_F$ the Frobenius norm of a matrix, defined by
\[
\|A\|_F=\left(\sum_{i,j=1}^m|a_{ij}|^2\right)^{1/2}.
\]
Two matrix norms satisfy $\|A\|\le\|A\|_F\le\sqrt m\,\|A\|$.

Let $\Omega$ be a bounded smooth domain and let
$\Omega_\varepsilon=\Phi_\varepsilon(\Omega)$, where
$\Phi_\varepsilon$ is a smooth family of diffeomorphisms satisfying
$\Phi_0(x)=x$. We set $V=\partial_\varepsilon\Phi_\varepsilon|_0$
and $h=V\cdot\nu$ on $\partial\Omega$. If we consider the radial deformation on the ball, we take
$\Phi_\varepsilon(\theta)=\rho_\varepsilon(\theta)\theta$ on $\partial B$ and $V(\theta)=h(\theta)\theta$.

Given a Dirichlet eigenvalue $\lambda_0$ of multiplicity $m$ and a real
$L^2(\Omega)$-orthonormal basis $u_1,\ldots,u_m$ of its eigenspace $E_0$, we denote by $E_\varepsilon$ the span of the
eigenfunctions whose eigenvalues converge to $\lambda_0$.
On the unit ball $B$, we take $\lambda_0 = \lambda_2(B)$ and write
$\mathcal N_*=\ker(-\Delta_B-\lambda_*)$ for $E_0$. The following lemma associates a matrix with the perturbed
eigenvalues and gives its second-order expansion.

\begin{lemma}\label{lem:orthogonal-cluster-basis}
There exists a $L^2(\Omega_\varepsilon)$-orthonormal basis
$U_1(\varepsilon),\ldots,U_m(\varepsilon)$ of $E_\varepsilon$ that
depends smoothly on $\varepsilon$ and satisfies $U_i(0)=u_i$.
The restriction of the Dirichlet Laplacian to $E_\varepsilon$ has a real
symmetric matrix $H_\varepsilon=(h_{ij}(\varepsilon))_{i,j=1}^m$ in this basis, determined by
\begin{equation}\label{eq:cluster-interior}
-\Delta U_i\left(\varepsilon\right)
=\sum_{j=1}^m h_{ji}\left(\varepsilon\right)U_j\left(\varepsilon\right)
\quad\text{in }\Omega_\varepsilon,
\qquad U_i\left(\varepsilon\right)=0
\quad\text{on }\partial\Omega_\varepsilon.
\end{equation}
Its eigenvalues are the perturbed eigenvalues.
Under the boundary
expansion~\eqref{eq:radialpath}, $U_i\left(\varepsilon\right)$ and $H_\varepsilon$ satisfy the following expansion
\begin{equation}\label{eq:cluster-taylor}
\begin{aligned}
U_i\left(\varepsilon,x\right)
 &=u_i\left(x\right)+\varepsilon v_i\left(x\right)
   +\varepsilon^2w_i\left(x\right)+o\left(\varepsilon^2\right),\\
H_\varepsilon
 &=\lambda_0 I_m+\varepsilon M\left(h\right)
   +\varepsilon^2 K_\Phi+o\left(\varepsilon^2\right).
\end{aligned}
\end{equation}
\end{lemma}

\begin{proof}
We give a brief proof. By perturbation theory for an isolated
spectral cluster~\cite{Kato}, we may choose an orthonormal basis
$U_1(\varepsilon),\ldots,U_m(\varepsilon)$ of $E_\varepsilon$
that depends smoothly on $\varepsilon$ and satisfies $U_i(0)=u_i$.
Since the Dirichlet Laplacian preserves \(E_\varepsilon\), we have $-\Delta E_\varepsilon\subseteq E_\varepsilon$. Hence, its restriction to this subspace admits a representation by a smooth real symmetric matrix \(H_\varepsilon\) in this basis, which yields \eqref{eq:cluster-interior}. Taylor expansion at $\varepsilon=0$
then gives~\eqref{eq:cluster-taylor}. 
\end{proof}

We next derive the interior equations and boundary conditions satisfied by \(v_i\) and \(w_i\). As in Henrot and Pierre~\cite[Section~5.6]{HP},
we differentiate the interior system at a fixed spatial point.
There the Laplacian commutes with differentiation
in $\varepsilon$. Substituting both lines of~\eqref{eq:cluster-taylor}
into~\eqref{eq:cluster-interior} and comparing the coefficients of
$\varepsilon$ and $\varepsilon^2$, respectively, gives
\begin{equation}\label{eq:cluster-coefficient-equations}
\begin{aligned}
\left(-\Delta-\lambda_0\right)v_i
 &=\sum_{j=1}^m m_{ji}\left(h\right)u_j,\\
\left(-\Delta-\lambda_0\right)w_i
 &=\sum_{j=1}^m m_{ji}\left(h\right)v_j+\sum_{j=1}^m k_{ji}u_j
\end{aligned}
\qquad\text{in }\Omega.
\end{equation}

For the boundary condition of $v_i$, we use the Dirichlet shape derivative
of Henrot and Pierre~\cite[Section~5.6, formula~(5.83)]{HP}.
Since the pullback $U_i(\varepsilon)\circ\Phi_\varepsilon$ vanishes
on $\partial\Omega$ for all $\varepsilon$, its derivative there is zero.
The chain rule therefore gives
\begin{equation}\label{eq:cluster-first-boundary}
v_i=-\nabla u_i\cdot V=-h\,\partial_\nu u_i
\qquad\text{on }\partial\Omega.
\end{equation}
Here the tangential gradient of $u_i$ vanishes since its trace is zero.
Multiplying the first equation in~\eqref{eq:cluster-coefficient-equations}
by $u_j$ and integrating over $\Omega$, we obtain the Hadamard formula
for multiple eigenvalues~\cite{ESI,HP,SZ} by Green's identity:
\begin{equation}\label{eq:hadamard}
m_{ji}\left(h\right)=\int_{\partial\Omega}v_i\,\partial_\nu u_j\,\mathrm{d}S
=-\int_{\partial\Omega}h\,\partial_\nu u_i\,\partial_\nu u_j\,\mathrm{d}S.
\end{equation}

Since the subsequent second-order calculations concern radial
deformations of the ball, we now restrict ourselves to
$\Omega=B$ and $\lambda_0=\lambda_*$ and derive the boundary
condition for $w_i$. For the boundary condition of $w_i$, fix
$\theta\in\mathbb S^{d-1}$ and set
$\delta_\varepsilon=\rho_\varepsilon(\theta)-1$.
By Lemma~\ref{lem:volume},
$\delta_\varepsilon=\varepsilon h+\varepsilon^2g+O(\varepsilon^3)$
and $\delta_\varepsilon^2=\varepsilon^2h^2+O(\varepsilon^3)$ at this $\theta$. By elliptic regularity up to the boundary and smooth
extension~\cite{GT,HP}, we may expand the functions in the radial
variable at $r=1$. For brevity, we omit $\theta$ from the notation.
All radial derivatives on the right are evaluated at $r=1$. Thus, we have
\[
\begin{aligned}
&u_i\left(1+\delta_\varepsilon\right)
 =u_i+\varepsilon h\,\partial_r u_i
   +\varepsilon^2\left(g\,\partial_r u_i+
              \tfrac12h^2\partial_{rr}u_i\right)+o\left(\varepsilon^2\right),\\
&\varepsilon v_i\left(1+\delta_\varepsilon\right)
 =\varepsilon v_i+\varepsilon^2h\,\partial_r v_i+o\left(\varepsilon^2\right),\\
&\varepsilon^2 w_i\left(1+\delta_\varepsilon\right)
 =\varepsilon^2w_i+o\left(\varepsilon^2\right).
\end{aligned}
\]
Since $U_i(\varepsilon,\rho_\varepsilon(\theta)\theta)=0$,
the coefficient of $\varepsilon^2$ in its Taylor expansion
must vanish. Hence we obtain
\begin{equation}\label{eq:cluster-second-boundary}
w_i=-h\,\partial_r v_i-g\,\partial_r u_i
    -\frac{h^2}{2}\partial_{rr}u_i
\qquad\text{on }\mathbb S^{d-1}.
\end{equation}

\begin{remark}
We call $M(h)=\frac{d}{d\varepsilon}H_\varepsilon|_{\varepsilon=0}$ \emph{the first-order matrix} and $K_\Phi=\tfrac12\frac{d^2}{d\varepsilon^2}H_\varepsilon|_{\varepsilon=0}$
\emph{the second-order matrix} along the deformation $\Phi_\varepsilon$ in the chosen basis. Here $M(h)=(m_{ij}(h))$ and $K_\Phi=(k_{ij})$. The eigenvalues of $M(h)$ are the first-order coefficients in the
expansions of the perturbed eigenvalues, after suitable labeling.
If $M(h)=0$, the eigenvalues of $K_\Phi$ are the corresponding
second-order coefficients.
\end{remark}

\subsection{The reduced boundary operator}\label{subsec:reduced}
We introduce a boundary operator for the subsequent calculations.
Since $\lambda_*$ is a Dirichlet eigenvalue of the unit ball, the boundary problem for $-\Delta-\lambda_*$ is solvable only for boundary data satisfying a compatibility condition. Therefore, we use a reduced Dirichlet-to-Neumann operator. Let
\[
\mathcal X=\left\{f\in H^{1/2}\left(\mathbb{S}^{d-1}\right):\Pi_1f=0\right\}.
\]
For $f\in\mathcal X$, let $v_f$ be the solution of
\[
\begin{cases}
	-\Delta v_f=\lambda_*v_f
	& \text{in }B,\\
	v_f=f
	& \text{on }\mathbb{S}^{d-1}
\end{cases}
\]
subject to the orthogonality condition
$$
\displaystyle\int_B v_fu\,\mathrm{d}x=0
\quad \text{for all }\quad u\in\mathcal N_*.
$$
The Fredholm alternative makes this solution unique. The required elliptic solvability and trace theory can be found in~\cite{McLean}. We define
\[
Tf=\partial_{\nu} v_f\in H^{-1/2}\left(\mathbb{S}^{d-1}\right).
\]
The operator is symmetric in the boundary duality because Green's formula gives
\[
\left\langle f,Tg\right\rangle
=\int_B\left(\nabla v_f\cdot\nabla v_g-\lambda_*v_fv_g\right)\,\mathrm{d} x
=\left\langle g,Tf\right\rangle.
\]

By direct calculation, for $Y\in\mathcal{H}_\ell$, $\ell\ne1$, separation of variables gives $v_Y(r\theta)=\varphi_\ell(r)Y(\theta)$, where
\[
\varphi_\ell\left(r\right)=r^{1-d/2}\frac{J_{\ell+d/2-1}\left(j_{d/2,1}r\right)}{J_{\ell+d/2-1}\left(j_{d/2,1}\right)}.
\]
By the monotonicity and interlacing properties of Bessel zeros,
the denominator is nonzero for every $\ell\ne1$.
We define
\[
t_\ell:=\varphi_\ell'(1),\qquad \ell\ne1.
\]
By the definition of $T$, for every $Y\in\mathcal{H}_\ell$
with $\ell\ne1$, we have $TY=\varphi_\ell'(1)Y=t_\ell Y$. Thus $T$ acts on $\mathcal{H}_\ell$ as multiplication by $t_\ell$.
For $\ell=0$, the Bessel derivative recurrence and
$J_{d/2}(j_{d/2,1})=0$ give $t_0=0$. The preceding representation allows us to compute \(t_\ell\) explicitly. For integers $n\ge3$, we also define the coefficients \(A_{d,n}\) that appear in the trace theorem below:
\begin{equation}\label{eq:coefficients}
t_\ell
=j_{d/2,1}\frac{J'_{\ell+d/2-1}\left(j_{d/2,1}\right)}{J_{\ell+d/2-1}\left(j_{d/2,1}\right)}-\left(\frac{d}{2}-1\right),\qquad
A_{d,n}
=d-1+\frac{nt_{n-1}+\left(n+d-2\right)t_{n+1}}{2n+d-2}.
\end{equation}
The estimates for the Bessel zero and the coefficients $t_\ell$ used below are collected in Lemmas~\ref{lem:bracket} and~\ref{lem:tcoeff} in Appendix~\ref{app:bessel}.

\section{Proof of Theorem~\ref{thm:kernel}}\label{sec:trace}
We now prove Theorem~\ref{thm:kernel}. We first identify the first-order matrix and its kernel.
\paragraph{The first-order matrix.}
We choose an $L^2(B)$-orthonormal basis $u_1,\ldots,u_d$ of $\mathcal N_*$. They can be chosen so that
\begin{equation}\label{eq:normalization}
\partial_{\nu} u_i=\gamma\theta_i,
\end{equation}
where $\gamma^2=\frac{2d\lambda_*}{s_d}$. Indeed, by the Rellich identity~\cite{PTZ}, we have $\int_{\mathbb{S}^{d-1}}
\left(\partial_{\nu}u_i\right)^2\,\mathrm{d}S
=2\lambda_*$. Since $\int_{\mathbb{S}^{d-1}}\theta_i\theta_j\,\mathrm{d} S=\delta_{ij}s_d/d$, it gives the stated value of $\gamma^2$. Here $\delta_{ij}$ denotes the Kronecker delta. We apply~\eqref{eq:hadamard} to a zero-mean radial deformation. The first-order matrix $M(h)$
has entries
\begin{equation}\label{eq:firstmatrix}
	m_{ij}\left(h\right)
	=-\gamma^2\int_{\mathbb{S}^{d-1}}
	h\theta_i\theta_j\,\mathrm{d}S.
\end{equation}

We denote the trace of a matrix by $\operatorname{tr}$. Since
\[
\theta_i\theta_j=\frac{\delta_{ij}}d+
\left(\theta_i\theta_j-\frac{\delta_{ij}}d\right),
\]
and the second term belongs to $\mathcal{H}_2$, the matrix $M(h)$ depends only on $h_2$. Moreover,
\begin{equation}\label{eq:zero-trace}
\operatorname{tr} M\left(h\right)=-\gamma^2\int_{\mathbb{S}^{d-1}}h\sum_{i=1}^d\theta_i^2\,\mathrm{d} S
=-\gamma^2\int_{\mathbb{S}^{d-1}}h\,\mathrm{d} S=0.
\end{equation}
Thus, if $h$ has zero mean, then $M(h)$ has zero trace. We next give a quantitative lower bound for the largest eigenvalue
of $M(h)$ and characterize all perturbations $h$ for which $M(h)=0$.

\begin{lemma}\label{lem:firstmatrix}
For every real-valued $z\in\mathcal{H}_2$, there is a unique real
symmetric traceless matrix $Z=(z_{ij})_{i,j=1}^d$ such that
$z(\theta)=\theta^{\top}Z\theta$, where the superscript
$\top$ denotes transpose. We have
\begin{equation}\label{eq:matrixZ}
	M\left(z\right)=-\frac{4\lambda_*}{d+2}Z,
	\qquad
	\left\lVert z\right\rVert_2^2
	=\frac{2s_d}{d\left(d+2\right)}
	\operatorname{tr}\left(Z^2\right).
\end{equation}
In particular, $M(h)=0$ for a zero-mean $h$ if and only if $h_2=0$.
There is a constant $m_d>0$ such that
\begin{equation}\label{eq:positivefirst}
	\lambda_{\max}\left(M\left(z\right)\right)
	\ge m_d\left\lVert z\right\rVert_{1/2,*}
	\qquad z\in\mathcal{H}_2,
\end{equation}
where $\lambda_{\max}$ denotes the largest eigenvalue of a
symmetric matrix.
\end{lemma}

\begin{proof}
Every real homogeneous quadratic polynomial can be uniquely
written as $x^{\top}Zx$, where $Z$ is a real symmetric matrix.
Since $\Delta(x^{\top}Zx)=2\operatorname{tr}Z$, the polynomial
is harmonic if and only if $\operatorname{tr}Z=0$.
By homogeneity, the polynomial is uniquely determined by its
restriction to the unit sphere. By~\cite{ABR}, we have the following integral identity on the sphere:
\begin{equation}\label{eq:sphere-fourth-moment}
	\int_{\mathbb{S}^{d-1}}
	\theta_i\theta_j\theta_k\theta_l\,\mathrm{d}S
	=\frac{s_d}{d\left(d+2\right)}
	\left(
	\delta_{ij}\delta_{kl}
	+\delta_{ik}\delta_{jl}
	+\delta_{il}\delta_{jk}
	\right).
\end{equation}
Since $z(\theta)=\theta^\top Z\theta
=\sum_{k,l=1}^d z_{kl}\theta_k\theta_l$, substituting this
representation into~\eqref{eq:firstmatrix}, we obtain
\begin{align*}
	m_{ij}\left(z\right)
	&=-\gamma^2\sum_{k,l=1}^d z_{kl}
	\int_{\mathbb{S}^{d-1}}
	\theta_i\theta_j\theta_k\theta_l\,\mathrm{d}S\\
	&=-\frac{\gamma^2s_d}{d\left(d+2\right)}
	\left(\delta_{ij}\operatorname{tr}Z+z_{ij}+z_{ji}\right)
	=-\frac{4\lambda_*}{d+2}z_{ij}.
\end{align*}
Here the second equality follows from~\eqref{eq:sphere-fourth-moment},
and the third uses the symmetry and zero trace of $Z$, together
with $\gamma^2=2d\lambda_*/s_d$. Similarly, expanding $z^2$ gives
\begin{align*}
	\left\lVert z\right\rVert_2^2
	&=\sum_{i,j,k,l=1}^d z_{ij}z_{kl}
	\int_{\mathbb{S}^{d-1}}
	\theta_i\theta_j\theta_k\theta_l\,\mathrm{d}S\\
	&=\frac{s_d}{d\left(d+2\right)}
	\left(
	\left(\operatorname{tr}Z\right)^2
	+2\operatorname{tr}\left(Z^2\right)
	\right)=\frac{2s_d}{d\left(d+2\right)}
	\operatorname{tr}\left(Z^2\right).
\end{align*}
This proves~\eqref{eq:matrixZ}.
For zero-mean $h$, we have $M(h)=M(h_2)$ because
$\theta_i\theta_j\in\mathcal{H}_0\oplus\mathcal{H}_2$.
The first identity in~\eqref{eq:matrixZ} therefore gives
$M(h)=0$ if and only if $h_2=0$.

We now prove the last assertion.
Let $S$ be a real symmetric traceless $d\times d$ matrix,
with eigenvalues $\mu_1,\ldots,\mu_d$, and put
$L=\lambda_{\max}(S)\ge0$.
Since $\mu_i=-\sum_{j\ne i}\mu_j$, we have
$-(d-1)L\le\mu_i\le L$. Hence its Frobenius norm satisfies $\left\lVert S\right\rVert_{\mathrm F}^2
=\operatorname{tr}\left(S^2\right)
=\sum_{i=1}^d\mu_i^2
\le d\left(d-1\right)^2L^2$. Since $\|z\|_{1/2,*}^2=3\|z\|_2^2$ for $z\in\mathcal{H}_2$, \eqref{eq:matrixZ} gives
\[
\left\lVert M\left(z\right)\right\rVert_{\mathrm F}
=\frac{4\lambda_*}{d+2}
\left(\frac{d\left(d+2\right)}{6s_d}\right)^{1/2}
\left\lVert z\right\rVert_{1/2,*}.
\]
Applying the preceding bound to $S=M(z)$
proves~\eqref{eq:positivefirst}, for example with $m_d=\frac{4\lambda_*}
{\left(d-1\right)\sqrt{6s_d\left(d+2\right)}}$.
\end{proof}

\begin{remark}
\eqref{eq:firstmatrix} defines $M(f)$ for every
$f\in L^2(\mathbb{S}^{d-1})$. By the Cauchy--Schwarz inequality, $|m_{ij}(f)|
\le \gamma^2\|f\|_2\|\theta_i\theta_j\|_2
\le C_d\|f\|_2$. Since $\|M(f)\|\le\|M(f)\|_F
\le d\max_{1\le i,j\le d}|m_{ij}(f)|$, we obtain
$\|M(f)\|\le C_d\|f\|_2$. In particular, by~\eqref{eq:hnorm},
\begin{equation}\label{eq:hadamard-norm}
	\|M(f)\|\le C_d\|f\|_2\le C_d\|f\|_{1/2,*}
	\qquad \text{for every } f\in H^{1/2}(\mathbb{S}^{d-1}).
\end{equation}
If $f\equiv c$ for some constant $c\in\mathbb{R}$, we obtain
\begin{equation}\label{constant-M(c)}
	M(c)=-\gamma^2c\,\frac{s_d}{d}I_d=-2\lambda_*cI_d.
\end{equation}
\end{remark}

\paragraph{The second-order matrix.}
For a volume-preserving radial deformation, $h_0=0$ by
Lemma~\ref{lem:volume}. If also $h_2=0$, then $M(h)=0$ by
Lemma~\ref{lem:firstmatrix}. We next compute the second-order matrix of
Section~\ref{subsec:perturbation} for these directions.

\begin{proposition}\label{prop:secondmatrix}
Consider a smooth volume-preserving radial deformation as in
Lemma~\ref{lem:volume}, with $h_2=0$. In the basis~\eqref{eq:normalization},
the second-order matrix $K_\Phi$ has entries
\begin{equation}\label{eq:secondmatrix}
k_{ij}=\gamma^2\left(
\left\langle h\theta_j,T\left(h\theta_i\right)\right\rangle
+\int_{\mathbb{S}^{d-1}}\left(\frac{d-1}{2}h^2-g\right)\theta_i\theta_j\,\mathrm{d} S\right).
\end{equation}
Its trace satisfies
\begin{equation}\label{eq:tracematrix}
\operatorname{tr} K_\Phi=\gamma^2\left(
\sum_{i=1}^d\left\langle h\theta_i,T\left(h\theta_i\right)\right\rangle
+\left(d-1\right)\left\lVert h\right\rVert_2^2\right).
\end{equation}
\end{proposition}
\begin{proof}
We first check compatibility of the boundary data. For every $j$, since $h_0=h_2=0$, we obtain $\int_{\mathbb{S}^{d-1}}h\theta_i\theta_j\,\mathrm{d} S=0$. Thus $\Pi_1(h\theta_i)=0$ and $h\theta_i\in\mathcal X$.

To simplify the proof, we may choose the basis $U_i(\varepsilon)$
in Lemma~\ref{lem:orthogonal-cluster-basis} so that
$v_i\perp\mathcal N_*$. Indeed, differentiating 
$$
\int_{\Omega_\varepsilon}U_iU_j\,\mathrm{d}x
=\delta_{ij}
$$
with respect to $\varepsilon$ at $\varepsilon=0$ and using
$u_i=u_j=0$ on $\partial B$, we obtain 
$$
\int_B\left(v_i u_j+u_i v_j\right)\,\mathrm{d}x=0.
$$
Hence the matrix
$\mathcal B=\left(\int_B u_i v_j\,\mathrm{d}x\right)_{i,j=1}^d$
is skew-symmetric. Then the basis $(U_1,\ldots,U_d)\exp(-\varepsilon\mathcal B)$ eliminates
the projections of $v_i$ onto $\mathcal N_*$.
Thus we may assume that $v_i\perp\mathcal N_*$ for every $i$. Since $M(h)=0$, the two interior
equations~\eqref{eq:cluster-coefficient-equations} reduce to
\begin{equation}\label{eq:ball-correction-equations}
\left(-\Delta-\lambda_*\right)v_i=0,
\qquad
\left(-\Delta-\lambda_*\right)w_i=\sum_{j=1}^d k_{ji}u_j
\qquad\text{in }B.
\end{equation}
By~\eqref{eq:cluster-first-boundary}, \eqref{eq:normalization} and $v_i\perp\mathcal N_*$, the definition of $T$ gives
\begin{equation}\label{eq:firstcorrection}
v_i|_{\mathbb{S}^{d-1}}=-\gamma h\theta_i,
\qquad
\partial_{\nu} v_i=-\gamma T\left(h\theta_i\right).
\end{equation}
We also need the second radial derivative of $u_i$. The interior equation in polar coordinates is
\[
-\partial_{rr}u_i-\frac{d-1}{r}\partial_ru_i
-\frac1{r^2}\Delta_{\mathbb{S}^{d-1}}u_i=\lambda_*u_i.
\]
At $r=1$, $u_i=0$, so we have $\Delta_{\mathbb{S}^{d-1}}u_i=0$. Thus $\partial_{rr}u_i|_{\mathbb{S}^{d-1}}
=-\left(d-1\right)\partial_ru_i|_{\mathbb{S}^{d-1}}
=-\left(d-1\right)\gamma\theta_i$. Combining this identity with~\eqref{eq:cluster-second-boundary}, \eqref{eq:normalization}, and~\eqref{eq:firstcorrection}, we obtain
\begin{equation}\label{eq:secondboundary}
w_i|_{\mathbb{S}^{d-1}}
=\gamma hT\left(h\theta_i\right)
+\gamma\left(\frac{d-1}{2}h^2-g\right)\theta_i.
\end{equation}
Multiplying the second equation in~\eqref{eq:ball-correction-equations}
by $u_l$ and integrating over $B$, we obtain
\begin{align*}
k_{li}
&=\int_Bu_l\left(-\Delta-\lambda_*\right)w_i\,\mathrm{d} x\\
&=\int_{\mathbb{S}^{d-1}}\left(w_i\partial_{\nu} u_l-u_l\partial_{\nu} w_i\right)\,\mathrm{d} S=\int_{\mathbb{S}^{d-1}}w_i\partial_{\nu} u_l\,\mathrm{d} S.
\end{align*}
By~\eqref{eq:normalization}, \eqref{eq:secondboundary}, and the symmetry of $T$, we obtain~\eqref{eq:secondmatrix}. We now sum the diagonal entries. Since $\sum_i\theta_i^2=1$, the second integral in the trace is
\[
\int_{\mathbb{S}^{d-1}}\left(\frac{d-1}{2}h^2-g\right)\,\mathrm{d} S
=\frac{d-1}{2}\left\lVert h\right\rVert_2^2-\int_{\mathbb{S}^{d-1}}g\,\mathrm{d} S
=\left(d-1\right)\left\lVert h\right\rVert_2^2
\]
by~\eqref{eq:volumeconstraint}. This proves~\eqref{eq:tracematrix}.
\end{proof}

The trace in~\eqref{eq:tracematrix} still contains the products
$h\theta_i$. We use the following lemma to eliminate the coordinate functions $\theta_i$ from~\eqref{eq:tracematrix} and express the trace in terms of the spherical harmonic components of $h$.

\begin{lemma}\label{lem:weights}
For $Y\in\mathcal{H}_n$, $n\ge1$, we have $\theta_iY\in\mathcal{H}_{n-1}\oplus\mathcal{H}_{n+1}$ and
\begin{equation}\label{eq:weights}
\sum_{i=1}^d\left\lVert \Pi_{n-1}\left(\theta_iY\right)\right\rVert_2^2
=\frac{n}{2n+d-2}\left\lVert Y\right\rVert_2^2,\qquad
\sum_{i=1}^d\left\lVert \Pi_{n+1}\left(\theta_iY\right)\right\rVert_2^2
=\frac{n+d-2}{2n+d-2}\left\lVert Y\right\rVert_2^2.
\end{equation}
If $h_0=h_1=h_2=0$, then
\begin{equation}\label{eq:diagonalform}
\sum_{i=1}^d\left\langle h\theta_i,T\left(h\theta_i\right)\right\rangle
=\sum_{n\ge3}
\frac{nt_{n-1}+\left(n+d-2\right)t_{n+1}}{2n+d-2}\left\lVert h_n\right\rVert_2^2.
\end{equation}
The identity holds for smooth $h$ and, by continuity, for every such $h\in H^{1/2}(\mathbb{S}^{d-1})$.
\end{lemma}
\begin{proof}
Let $P_n$ be the homogeneous harmonic polynomial of degree $n$
whose restriction to $\mathbb{S}^{d-1}$ is $Y$. We use the standard harmonic decomposition of a polynomial~\cite[Theorem~5.21]{ABR} in the particular form
\[
x_iP_n=H_{n+1,i}+\frac{\left|x\right|^2}{2n+d-2}\partial_iP_n,
\]
where $H_{n+1,i}$ is homogeneous and harmonic of degree $n+1$. The derivative $\partial_iP_n$ is harmonic of degree $n-1$. On the sphere this identifies
\[
\Pi_{n-1}\left(\theta_iY\right)
=\frac{\partial_iP_n|_{\mathbb{S}^{d-1}}}{2n+d-2}.
\]
By homogeneity, we have $P_n(r\theta)=r^nY(\theta)$,
which gives $\partial_rP_n=nY$ at $r=1$. Therefore, we have
\begin{align*}
\sum_{i=1}^d\left\lVert \partial_iP_n\right\rVert_2^2=\left\lVert\nabla P_n\right\rVert_2^2
&=n^2\left\lVert Y\right\rVert_2^2+\left\lVert \nabla_{\mathbb{S}^{d-1}}Y\right\rVert_2^2\\
&=\left[n^2+n\left(n+d-2\right)\right]\left\lVert Y\right\rVert_2^2\\
&=n\left(2n+d-2\right)\left\lVert Y\right\rVert_2^2.
\end{align*}
In the second line we used integration by parts on the sphere. Dividing by $(2n+d-2)^2$ proves the first formula in~\eqref{eq:weights}. Orthogonality and $\sum_i\theta_i^2=1$ imply
\begin{align*}
\sum_i\left\lVert \Pi_{n+1}\left(\theta_iY\right)\right\rVert_2^2
&=\sum_i\left\lVert \theta_iY\right\rVert_2^2
-\sum_i\left\lVert \Pi_{n-1}\left(\theta_iY\right)\right\rVert_2^2\\
&=\left(1-\frac{n}{2n+d-2}\right)\left\lVert Y\right\rVert_2^2,
\end{align*}
which is the second formula.

Multiplication by \(\theta_i\) maps a spherical harmonic of degree $n$ into \(\mathcal H_{n-1}\oplus\mathcal H_{n+1}\). Therefore, if $h$ is expanded in spherical harmonics, cross terms between harmonics of degree $n$ and \(n+2\) may appear: both products contain a component of degree \(n+1\). We now show that the sum of these cross terms over $i$ is zero. To see this, take \(Y\in\mathcal H_n\) and \(Z\in\mathcal H_{n+2}\). Let $P_{n+2}$ be the homogeneous harmonic polynomial of degree ${n+2}$
whose restriction to $\mathbb S^{d-1}$ is $Z$. The only common degree is $n+1$, and
$$
\sum_i\left\langle \Pi_{n+1}\left(\theta_iY\right),\Pi_{n+1}\left(\theta_iZ\right)\right\rangle=\frac1{2n+d+2}
\sum_i\int_{\mathbb{S}^{d-1}}
\left(\theta_iY-\frac{\partial_iP_n}{2n+d-2}\right)
\partial_iP_{n+2}\,\mathrm{d} S.
$$
Euler's identity for homogeneous polynomials gives
\[
\sum_i\int_{\mathbb{S}^{d-1}}\theta_iY\partial_iP_{n+2}\,\mathrm{d} S
=\left(n+2\right)\int_{\mathbb{S}^{d-1}}YZ\,\mathrm{d} S=0.
\]
For the other term, we split the gradients into radial and tangential components:
\begin{align*}
\sum_i\int_{\mathbb{S}^{d-1}}\partial_iP_n\partial_iP_{n+2}\,\mathrm{d} S
&=n\left(n+2\right)\int_{\mathbb{S}^{d-1}}YZ\,\mathrm{d} S
+\int_{\mathbb{S}^{d-1}}\nabla_{\mathbb{S}^{d-1}}Y\cdot\nabla_{\mathbb{S}^{d-1}}Z\,\mathrm{d} S\\
&=\left[n\left(n+2\right)+\left(n+2\right)\left(n+d\right)\right]
\int_{\mathbb{S}^{d-1}}YZ\,\mathrm{d} S=0.
\end{align*}
Thus there are no cross terms in the sum over $i$.

Let $h=\sum_{n=3}^N h_n$ be a finite spherical harmonic sum.
For each $n$, the product $h_n\theta_i$ has components of degrees
$n-1$ and $n+1$. These components are orthogonal, and $T$ acts on
them by multiplication by $t_{n-1}$ and $t_{n+1}$, respectively.
Thus, by~\eqref{eq:weights},
\begin{align*}
	\sum_{i=1}^d \langle h_n\theta_i,T(h_n\theta_i)\rangle
	&=t_{n-1}\sum_{i=1}^d
	\bigl\|\Pi_{n-1}(h_n\theta_i)\bigr\|_2^2
	+t_{n+1}\sum_{i=1}^d
	\bigl\|\Pi_{n+1}(h_n\theta_i)\bigr\|_2^2\\
	&=\frac{n t_{n-1}+(n+d-2)t_{n+1}}{2n+d-2}\|h_n\|_2^2.
\end{align*}
Since the cross terms vanish after summation over $i$, summing this
identity over $n$ gives~\eqref{eq:diagonalform} for $h$. Since both sides of~\eqref{eq:diagonalform} are continuous in the
$H^{1/2}$ norm, the identity extends by density from finite spherical
harmonic sums to every $h\in H^{1/2}(\mathbb S^{d-1})$ with $h_0=h_1=h_2=0$.
\end{proof}

\begin{proof}[Proof of Theorem~\ref{thm:kernel}]
Along the deformation $\Omega_\varepsilon=\Omega_{\varepsilon h}$ of~\eqref{eq:radialpath}, Proposition~\ref{prop:secondmatrix} and Lemma~\ref{lem:weights} imply
\begin{align*}
\frac1d\operatorname{tr} K_\Phi
&=\frac{\gamma^2}{d}\sum_{n\ge3}
\left(d-1+\frac{nt_{n-1}+\left(n+d-2\right)t_{n+1}}{2n+d-2}\right)
\left\lVert h_n\right\rVert_2^2\\
&=\frac{2\lambda_*}{s_d}\sum_{n\ge3}A_{d,n}\left\lVert h_n\right\rVert_2^2.
\end{align*}
By Lemma~\ref{lem:orthogonal-cluster-basis}, $\sum_{j=2}^{d+1}\lambda_j(\Omega_\varepsilon)
=\operatorname{tr}H_\varepsilon$. Taking the trace in~\eqref{eq:cluster-taylor} and using \eqref{eq:zero-trace}, we obtain
\[
\frac{1}{d}\sum_{j=2}^{d+1}\lambda_j(\Omega_\varepsilon)
=\lambda_*+\frac{\varepsilon^2}{d}\operatorname{tr}K_\Phi
+o(\varepsilon^2).
\]
Thus we obtain~\eqref{eq:maintrace}.

We compute the lowest coefficient separately. Let $\Delta_d=\left(d+4\right)\left(d+2\right)-\lambda_*>0$. A direct calculation using Lemma~\ref{lem:bracket} and~\eqref{eq:lowcoeff} gives
$$
A_{d,3}
=\frac{\left(d+2\right)\left(\lambda_*-3\left(d+2\right)\right)}{\Delta_d}>0.
$$
For $n\ge4$, both indices $n-1$ and $n+1$ are at least three. By Lemma~\ref{lem:tcoeff},
\[
A_{d,n}\ge d-1+t_3
=\frac{\lambda_*}{d+2}-2>1.
\]
We also use~\eqref{eq:tasymptotic} to calculate the growth. The leading weighted average is
\[
\frac{n\left(n-1\right)+\left(n+d-2\right)\left(n+1\right)}{2n+d-2}
=n+\frac{d-2}{2n+d-2}.
\]
The weighted errors are $O_d(n^{-1})$. Hence $A_{d,n}=n+d-1+O_d\left(n^{-1}\right)$. There exists an integer $N_d\ge4$ for which $A_{d,n}/(1+n)\ge1/2$ whenever $n\ge N_d$. We may therefore choose
\[
c_d=\min\left\{\frac12,\frac{A_{d,3}}4,
\frac{A_{d,4}}5,\ldots,\frac{A_{d,N_d-1}}{N_d}\right\}>0.
\]
This proves~\eqref{eq:coercivity}. Finally, we use
\[
\lambda_{d+1}\left(\Omega_{\varepsilon h}\right)
\ge\frac1d\sum_{j=2}^{d+1}\lambda_j\left(\Omega_{\varepsilon h}\right).
\]
Formula~\eqref{eq:maintrace} and the bound~\eqref{eq:coercivity} prove the last assertion.
\end{proof}

\section{Proof of Theorem~\ref{thm:uniform} and its corollary}\label{sec:uniform}
In Section~\ref{sec:trace}, we give
the second-order expansion of the mean eigenvalue of the cluster
along each fixed smooth deformation direction $h$ with $h_0=h_1=h_2=0$. We now establish an estimate valid in one common neighborhood. We first establish a lemma that provides uniform estimates for general deformations of an arbitrary smooth domain. For a vector field $V=(V_1,\ldots,V_d)$, its Jacobian
matrix is $DV=(\partial_k V_i)_{i,k=1}^d$. We use the Frobenius norm
pointwise when taking norms of matrix fields. In particular,
\[
\lVert DV\rVert_{L^\infty(\Omega)}
=\operatorname*{ess\,sup}_{x\in\Omega}\lVert DV(x)\rVert_{\mathrm F},
\qquad
\lVert DV\rVert_{L^2(\Omega)}
=\left(\int_\Omega \lVert DV(x)\rVert_{\mathrm F}^2\,\mathrm{d}x\right)^{1/2}.
\]

\begin{lemma}\label{lem:cluster}
Let $\Omega\subset\mathbb{R}^d$ be a bounded smooth domain, and let $\lambda_0>0$ be an isolated Dirichlet eigenvalue of multiplicity $m$. Suppose that $V\in C^1(\overline{\Omega};\mathbb{R}^d)$ and the map $F(x)=x+V(x)$ is a diffeomorphism onto its image. We set
\[
\delta=\left\lVert DV\right\rVert_{L^\infty\left(\Omega\right)},
\qquad \sigma=\left\lVert DV\right\rVert_{L^2\left(\Omega\right)}.
\]
For sufficiently small $\delta$, the eigenvalues of the cluster near $\lambda_0$ on $F(\Omega)$ are the eigenvalues of a real symmetric $m\times m$ matrix
\begin{equation}\label{eq:clusterexpansion}
H_\Omega\left(V\right)
=\lambda_0I_m+H_{\Omega,1}\left(V\right)
+H_{\Omega,2}\left(V,V\right)+\mathcal R_\Omega\left(V\right).
\end{equation}
The map $H_{\Omega,1}$ is linear in $V$, and $H_{\Omega,2}$ is a symmetric bilinear map. Both maps take values in the real symmetric $m\times m$ matrices. There are constants independent of $V$ such that
\begin{equation}\label{eq:clusterbounds}
\left\lVert H_{\Omega,2}\left(V,V\right)\right\rVert\le C\sigma^2,
\qquad
\left\lVert \mathcal R_\Omega\left(V\right)\right\rVert\le C\delta\sigma^2.
\end{equation}
\end{lemma}
\begin{proof}

\noindent\textbf{Step 1. Reformulation on the fixed domain.}
We take \(\mathcal E = H_0^1(\Omega)\) to be the Hilbert space equipped with the energy inner product
\((u,v)_{\mathcal E} = \int_\Omega \nabla u\cdot\nabla v\,\mathrm{d}x.\)
By Poincar\'e's inequality, the associated norm \(\|u\|_{\mathcal E} = \bigl(\int_\Omega |\nabla u|^2\,\mathrm{d}x\bigr)^{1/2}\) is equivalent to the usual \(H^1_0(\Omega)\) norm. Let $G=DV$ and define the scalar and matrix coefficients
\[
\mathcal J(G)=\det(I_d+G),
\qquad
S(G)=\mathcal J(G)(I_d+G)^{-1}(I_d+G)^{-\top}.
\]
If $u=\widetilde u\circ F$, the chain rule and the change of variables
$y=F(x)$ give
\begin{align*}
\int_{F(\Omega)} |\nabla\widetilde u|^2\,\mathrm{d}y
=\int_\Omega (S(G)\nabla u)\cdot\nabla u\,\mathrm{d}x,\qquad
\int_{F(\Omega)} |\widetilde u|^2\,\mathrm{d}y
=\int_\Omega \mathcal J(G)u^2\,\mathrm{d}x.
\end{align*}
By the Riesz representation theorem, the corresponding bilinear forms
define self-adjoint operators $A(V)$ and $B(V)$ on
$\mathcal E$ by the identities
\begin{align}
(A(V)u,v)_{\mathcal E}
&=\int_\Omega (S(G)\nabla u)\cdot\nabla v\,\mathrm{d}x,
\label{eq:stiffness-form}\\
(B(V)u,v)_{\mathcal E}
&=\int_\Omega \mathcal J(G)uv\,\mathrm{d}x.
\label{eq:mass-form}
\end{align}
For sufficiently small $\delta$, the bilinear form associated with
$A(V)$ is coercive. Thus $A(V)$ is positive, and the Lax--Milgram theorem yields
a bounded inverse on $\mathcal E$. The operator $B(V)$ is compact because the embedding
$H_0^1(\Omega)\hookrightarrow L^2(\Omega)$ is compact. We have $A(0)=\operatorname{Id}_{\mathcal E}$. We denote $B(0)$ by $B_0$. The weak form of the eigenvalue equation is
\[
A\left(V\right)u=\lambda B\left(V\right)u.
\]
For $w=A(V)^{1/2}u$, this becomes
\[
T_{\Omega}\left(V\right)w=\lambda^{-1}w,
\]
where $T_{\Omega}\left(V\right)=A\left(V\right)^{-1/2}B\left(V\right)A\left(V\right)^{-1/2}$.
The operator $T_{\Omega}(V)$ is compact and self-adjoint. Thus its nonzero eigenvalues are exactly the reciprocals of the Dirichlet eigenvalues on $F(\Omega)$.

\medskip
\noindent\textbf{Step 2. Some estimates.}
Fix $V$ and replace it by $tV$, so that $D(tV)=tG$.
Expanding the determinant as a polynomial in $t$ gives
\begin{equation}\label{J-expend}
\mathcal J(tG)=\det(I_d+tG)=1+\sum_{j=1}^{d}t^j\mathcal J_j(G).
\end{equation}
For $|t|\lVert G\rVert_{L^\infty}<1$, we have
$$
(I_d+tG)^{-1}=\sum_{p=0}^{\infty}(-t)^pG^p.
$$
Thus we may write
\begin{equation}\label{S-expend}
S(tG)=I_d+\sum_{j=1}^{\infty}t^jS_j(G).
\end{equation}
The $\mathcal J_j(G)$ and $S_j(G)$ are the coefficients
of $t^j$ in the respective expansions. They are homogeneous of degree $j$ in $G$. In particular, by the determinant expansion formula, we have
\begin{align*}
\mathcal J_1(G)=\operatorname{tr}G, \qquad
\mathcal J_2(G)=\tfrac12\bigl((\operatorname{tr}G)^2-\operatorname{tr}(G^2)\bigr).
\end{align*}
Moreover, we obtain
\begin{align*}
S_1(G)&=\mathcal J_1(G)I_d-G-G^{\top},\\
S_2(G)&=\mathcal J_2(G)I_d-\mathcal J_1(G)(G+G^{\top})
       +G^2+(G^{\top})^2+GG^{\top}.
\end{align*}
We set $\mathcal J_j=0$ for $j>d$. In \eqref{S-expend},
a term of degree $j$ comes from degrees $r,p,q$ with $r+p+q=j$
and $0\le r\le d$. For each $r$ there are at most $j+1$ choices of
$(p,q)$. By homogeneity and compactness of the unit sphere in
$\mathbb R^{d\times d}$, we have $|\mathcal J_r(G)|\le C_d\lVert G\rVert_{\mathrm F}^r, 1\le r\le d$. Since the determinant coefficients are homogeneous polynomials
and there are only finitely many terms, we have
\[
\lVert S_j(G)\rVert_{\mathrm F}+|\mathcal J_j(G)|
\le C_d(j+1)\lVert G\rVert_{\mathrm F}^j.
\]
Using $j+1\le 2^j$ for $j\ge1$, we obtain
\[
\lVert S_j(G)\rVert_{\mathrm F}+|\mathcal J_j(G)|
\le C^j\lVert G\rVert_{\mathrm F}^j
\qquad j\ge1,
\]
where $C$ depends only on $d$.

For fixed $u,v\in\mathcal E$, substitute \eqref{J-expend} and \eqref{S-expend} into
\eqref{eq:stiffness-form} and \eqref{eq:mass-form} with $V$ replaced by $tV$.
Collecting the coefficients of \(t^j\), we define $A_j(V)$ and
$B_j(V)$ by
\begin{equation}\label{eq:form-coefficients}
	(A_j(V)u,v)_{\mathcal E}
	=\int_\Omega (S_j(G)\nabla u)\cdot\nabla v\,\mathrm{d}x,\qquad
	(B_j(V)u,v)_{\mathcal E}
	=\int_\Omega \mathcal J_j(G)uv\,\mathrm{d}x.
\end{equation}
Moreover, we have
$$
A(tV)=\operatorname{Id}_{\mathcal E}
+\sum_{j\ge1}t^jA_j(V),\qquad
B(tV)=B_0+\sum_{j=1}^{d}t^jB_j(V).
$$
Since $\delta=\lVert G\rVert_{L^\infty}$ is sufficiently small, we may
set $t=1$ in these expansions. We thus obtain
\[
A\left(V\right)=\operatorname{Id}_{\mathcal E}+\sum_{j\ge1}A_j\left(V\right),
\qquad
B\left(V\right)=B_0+\sum_{j\ge1}B_j\left(V\right),
\]
where $A_j$ and $B_j$ are defined by \eqref{eq:form-coefficients}.

We let $\mu_0=1/\lambda_0$ and $P_0$ be the orthogonal projection in $\mathcal E$ onto the eigenspace of $B_0$ at $\mu_0$. Set $E_0=P_0\mathcal E$. Since the domain is smooth, elliptic regularity~\cite{GT} implies that the fixed eigenfunctions and their gradients are bounded on $\overline{\Omega}$. The eigenspace is finite-dimensional, so there is a constant $C_{E_0}$ such that
\begin{equation}\label{eq:fixedsmooth}
\left\lVert e\right\rVert_{L^\infty\left(\Omega\right)}+\left\lVert \nabla e\right\rVert_{L^\infty\left(\Omega\right)}
\le C_{E_0}\left\lVert e\right\rVert_{\mathcal E}
\qquad e\in P_0\mathcal E.
\end{equation}
We claim that
\begin{equation}
	\left\lVert A_j\right\rVert+\left\lVert B_j\right\rVert
	\le C^j\delta^j,
	\label{eq:operatorbound}
\end{equation}
\begin{equation}
	\left\lVert A_jP_0\right\rVert+\left\lVert B_jP_0\right\rVert
	\le C^j\delta^{j-1}\sigma,
	\label{eq:oneend}
\end{equation}
\begin{equation}
	\left\lVert P_0A_jP_0\right\rVert+\left\lVert P_0B_jP_0\right\rVert
	\le C^j\delta^{j-2}\sigma^2\qquad j\ge2.
	\label{eq:twoends}
\end{equation}
For the first estimate, since $\lVert S_j(G)\rVert_{\mathrm F}\le C^j\delta^j$ almost everywhere in $\Omega$, the Cauchy--Schwarz inequality gives $\left|\left(A_ju,v\right)_{\mathcal E}\right|
\le C^j\delta^j\left\lVert \nabla u\right\rVert_{L^2}\left\lVert \nabla v\right\rVert_{L^2}$. The same bound holds for $B_j$ by Poincar\'e's inequality. Using the definition of the norm and taking the supremum over all unit vectors \(u,v\in\mathcal E\) gives the desired operator norm bound. For the second estimate, let $j\ge 1$, $e\in P_0\mathcal E$, and
$v\in\mathcal E$. Then
\begin{align*}
	|(A_j e,v)_{\mathcal E}|
	=\left|\int_\Omega (S_j(G)\nabla e)\cdot\nabla v\,\mathrm{d}x\right|
	&\le \int_\Omega \|S_j(G)\|_{\mathrm F}\,|\nabla e|\,|\nabla v|\,\mathrm{d}x\\
	&\le C^j\|\nabla e\|_{L^\infty}
	\int_\Omega \|G\|_{\mathrm F}^j\,|\nabla v|\,\mathrm{d}x\\
	&\le C^j\|\nabla e\|_{L^\infty}
	\left(\int_\Omega \|G\|_{\mathrm F}^{2j}\,\mathrm{d}x\right)^{1/2}
	\|\nabla v\|_{L^2}\\
	&\le C^j\delta^{j-1}\sigma\,\|e\|_{\mathcal E}\|v\|_{\mathcal E}.
\end{align*}
The first line uses the inequality $|(A\xi)\cdot\eta|\le\|A\|_{\mathrm F}|\xi||\eta|$. For the fourth line, $\|G(x)\|_{\mathrm F}\le\delta$
gives
\[
\left(\int_\Omega \|G\|_{\mathrm F}^{2j}\,\mathrm{d}x\right)^{1/2}
\le \delta^{j-1}
\left(\int_\Omega \|G\|_{\mathrm F}^2\,\mathrm{d}x\right)^{1/2}
= \delta^{j-1}\sigma.
\]
The estimate for $B_j$ is similar. For the third estimate, let $j\ge 2$ and $u,v\in\mathcal E$.
Set $e=P_0u$ and $f=P_0v$. Then
\begin{align*}
	|(P_0A_jP_0u,v)_{\mathcal E}|
	=|(A_je,f)_{\mathcal E}|
	&\le \int_\Omega \|S_j(G)\|_{\mathrm F}\,|\nabla e|\,|\nabla f|
	\,\mathrm{d}x\\
	&\le C^j\|\nabla e\|_{L^\infty}\|\nabla f\|_{L^\infty}
	\int_\Omega \|G\|_{\mathrm F}^j\,\mathrm{d}x\\
	&\le C^j\|\nabla e\|_{L^\infty}\|\nabla f\|_{L^\infty}
	\delta^{j-2}\sigma^2\\
	&\le C^j\delta^{j-2}\sigma^2\|u\|_{\mathcal E}\|v\|_{\mathcal E}.
\end{align*}
The first line uses the self-adjointness of $P_0$. The third line uses
$\|G\|_{\mathrm F}^j\le\delta^{j-2}\|G\|_{\mathrm F}^2$ and
$\int_\Omega\|G\|_{\mathrm F}^2\,\mathrm{d}x=\sigma^2$. The estimate for $B_j$ is again similar.

\medskip
\noindent\textbf{Step 3. Matrix representation.}
In this step, the coefficients of \(t^j\) satisfy homogeneity of degree $j$ in $V$. For simplicity, we omit $V$ and, for example, we write \(T_{\Omega,j}\) for \(T_{\Omega,j}(V)\). We first expand
\[
A(tV)^{-1/2}
=\operatorname{Id}_{\mathcal E}+\sum_{j\ge1}t^j \widetilde{A}_j.
\]
The binomial series and the expansion of $A(tV)$ give
\[
\widetilde{A}_1=-\frac{1}{2}A_1,\qquad
\widetilde{A}_2=-\frac12A_2+\frac{3}{8}A_1^2.
\]
Substituting these formulas into
$T_{\Omega}(tV)=A(tV)^{-1/2}B(tV)A(tV)^{-1/2}$
and combining equal powers of $t$, we obtain the expansion of $T_{\Omega}(tV)$:
\[
T_{\Omega}(tV)=B_0+\sum_{j\ge1}t^jT_{\Omega,j},
\]
where
\begin{align*}
	T_{\Omega,1}
	&=B_1+\widetilde{A}_1B_0+B_0\widetilde{A}_1,\\
	T_{\Omega,2}
	&=B_2+\widetilde{A}_2B_0+B_0\widetilde{A}_2
	+\widetilde{A}_1B_1+B_1\widetilde{A}_1+\widetilde{A}_1B_0\widetilde{A}_1.
\end{align*}
We next construct a spectral projection whose range is the eigenspace associated with the perturbed eigenvalues. Let $\Gamma$ is a positively oriented circle around $\mu_0$
whose interior contains no other point of the spectrum of
$B_0$. Let
\[
R_0(z)=(z\operatorname{Id}_{\mathcal E}-B_0)^{-1},
\qquad D(t)=T_{\Omega}(tV)-B_0.
\]
Since $R_0(z)$ is uniformly bounded on $\Gamma$ and
$\|D(t)\|\le C|t|\delta$ for $|t|\delta$ small, we obtain
\[
(z\operatorname{Id}_{\mathcal E}-T_{\Omega}(tV))^{-1}
=\bigl(\operatorname{Id}_{\mathcal E}-R_0(z)D(t)\bigr)^{-1}R_0(z)
=\sum_{n\ge0}\bigl(R_0(z)D(t)\bigr)^nR_0(z)
\]
converges uniformly on $\Gamma$. Therefore, the Riesz projection
formula~\cite{Kato} gives
\[
P(tV)=\frac{1}{2\pi\mathrm i}\int_{\Gamma}
(z\operatorname{Id}_{\mathcal E}-T_{\Omega}(tV))^{-1}\,\mathrm dz
=P_0+\sum_{j\ge1}t^jP_j.
\]
Integrating the coefficients of $t$ and $t^2$, respectively, gives
\begin{align*}
	P_1
	&=\frac{1}{2\pi\mathrm i}\int_{\Gamma}
	R_0T_{\Omega,1}R_0\,\mathrm dz,\\
	P_2
	&=\frac{1}{2\pi\mathrm i}\int_{\Gamma}
	\bigl(R_0T_{\Omega,2}R_0+
	R_0T_{\Omega,1}R_0T_{\Omega,1}R_0\bigr)\,\mathrm dz,
\end{align*}
where $R_0=R_0(z)$. Comparing coefficients in $P(tV)^2=P(tV)$
and multiplying on both sides by $P_0$ yields $P_0P_1P_0=0, P_0P_2P_0=-P_0P_1^2P_0$. Now, let $E_t=P(tV)\mathcal E$ and $e_1,\ldots,e_m$ be a fixed orthonormal basis of $E_0$. For $|t|\delta$ sufficiently small, we have
$\|P(tV)-P_0\|<1$. Since $P_0e=e$ for every $e\in E_0$, $\|P(tV)e\|_{\mathcal E}
\ge \left(1-\|P(tV)-P_0\|\right)\|e\|_{\mathcal E}$. Thus $P(tV)|_{E_0}$ is injective. Since
$\dim E_0=\dim E_t=m$, it is an isomorphism
from $E_0$ onto $E_t$. Hence $P(tV)e_1,\ldots,P(tV)e_m$ form a basis of $E_t$. Note that the projected basis is not necessarily orthonormal. Therefore we introduce the Gram operator on $E_0$ by
\[
N(tV):=\left.P_0P(tV)P_0\right|_{E_0}.
\]
Indeed, since $P(tV)$ is an orthogonal projection and $P_0e_i=e_i$, we obtain $(N(tV)e_i,e_j)_{\mathcal E}
=(P(tV)e_i,P(tV)e_j)_{\mathcal E}$. Substituting the expansion of $P(tV)$ gives
\[
N(tV)=\operatorname{Id}_{E_0}
+\sum_{j\ge1}t^jN_j,
\qquad
N_j=\left.P_0P_jP_0\right|_{E_0}.
\]
Since $P_0P_1P_0=0$, we have $N_1=0$. Consequently,
\[
N(tV)^{-1/2}
=\operatorname{Id}_{E_0}-\frac12t^2N_2+O(t^3).
\]
Therefore, we define and expand
\begin{equation}
	W(tV):=\bigl(P(tV)|_{E_0}\bigr)N(tV)^{-1/2}
	=\left.\operatorname{Id}_{\mathcal E}\right|_{E_0}+\sum_{j\ge1}t^jW_j.
\end{equation}
Then $W(tV)$ maps $E_0$ onto $E_t$ and $W(tV)^*W(tV)
=N(tV)^{-1/2}N(tV)N(tV)^{-1/2}
=\operatorname{Id}_{E_0}$. Thus $W(tV)e_1,\ldots,W(tV)e_m$ form an orthonormal basis of $E_t$.
Multiplying the expansions of $P(tV)|_{E_0}$ and
$N(tV)^{-1/2}$ gives
$$
W_1=P_1|_{E_0}.
$$
In particular, $P_0W_1=(P_0P_1P_0)|_{E_0}=0$. Hence $\operatorname{ran} W_1\subseteq \ker P_0=E_0^\perp$.

Using the isometric isomorphism $W(tV)$, we can transfer
the restriction of $T_{\Omega}(tV)$ to $E_t$ onto the fixed
finite-dimensional space $E_0$, while preserving its
self-adjointness and eigenvalues. Therefore, we define and expand
\[
L(tV):=W(tV)^*T_{\Omega}(tV)W(tV)
=\mu_0\operatorname{Id}_{E_0}+\sum_{j\ge1}t^jL_j.
\]
Since $W(tV)$ maps $E_0$ onto $E_t$ and $T_{\Omega}(tV)$ leaves
$E_t$ invariant, we have $L(tV):E_0\longrightarrow E_0$. Its Taylor coefficients $L_j$ therefore map $E_0$ into itself,
so $P_0 L_j=L_j$ for every $j\ge 1$. In addition, it is easy to see that $T_{\Omega}(tV)W(tV)=W(tV)L(tV)$. Applying this identity to $e\in E_0$ and comparing the coefficients
of $t$ and $t^2$, we obtain
\begin{align*}
	B_0 W_1 e + T_{\Omega,1} e
	&= \mu_0 W_1 e + L_1 e,\\
	B_0 W_2 e + T_{\Omega,1} W_1 e + T_{\Omega,2} e
	&= \mu_0 W_2 e + W_1 L_1 e + L_2 e.
\end{align*}
Applying $P_0$ to both equations and using
$P_0B_0=\mu_0P_0$, $P_0 L_j=L_j$ and $P_0W_1=0$, we obtain
\[
L_1=P_0T_{\Omega,1}|_{E_0},\qquad
L_2=P_0T_{\Omega,2}|_{E_0}+P_0T_{\Omega,1}P_1|_{E_0}.
\]

Finally, since $W(tV)$ is an isometric isomorphism from $E_0$ onto
$E_t$, the operators $L(tV)$ and $T_\Omega(tV)|_{E_t}$
have the same eigenvalues. Hence the eigenvalues of
$L(tV)$ also are the reciprocals of the Dirichlet
eigenvalues in the cluster. We therefore get
\[
H_\Omega(tV)=L(tV)^{-1}
=\lambda_0\operatorname{Id}_{E_0}
+tH_{\Omega,1}(V)+t^2H_{\Omega,2}(V,V)+O(t^3).
\]
In a real orthonormal basis of $E_0$, this is a real symmetric
matrix with the perturbed eigenvalues. Comparing coefficients in
$L(tV)H_\Omega(tV)=\operatorname{Id}_{E_0}$, we get
\begin{equation}
	H_{\Omega,1}(V)=-\lambda_0^2L_1,\qquad
	H_{\Omega,2}(V,V)=\lambda_0^3L_1^2-\lambda_0^2L_2.
\end{equation}

For later use, we now compute $H_{\Omega,1}(V)$ explicitly. We choose a real $L^2(\Omega)$-orthonormal basis
$u_1,\ldots,u_m$ of the eigenspace at $\lambda_0$.
Setting $e_i=u_i/\sqrt{\lambda_0}$, we have
\[
(e_i,e_j)_{\mathcal E}
=\int_\Omega u_i u_j\,\mathrm{d}x=\delta_{ij}.
\]
Thus $e_1,\ldots,e_m$ form an orthonormal basis of $E_0$. Since
$B_0$ is self-adjoint and
$B_0 e_i=\lambda_0^{-1}e_i$, we have
\[
(T_{\Omega,1}e_i,e_j)_{\mathcal E}
=(B_1e_i,e_j)_{\mathcal E}
-\frac12\bigl((A_1B_0+B_0A_1)e_i,e_j\bigr)_{\mathcal E}
=(B_1e_i,e_j)_{\mathcal E}
-\frac{1}{\lambda_0}(A_1e_i,e_j)_{\mathcal E}.
\]
Using $H_{\Omega,1}(V)=-\lambda_0^2L_1$,
$L_1=P_0T_{\Omega,1}|_{E_0}$ and $P_0e_j=e_j$, we obtain
\begin{align*}
	(H_{\Omega,1}(V))_{ij}
	&= -\lambda_0^2(T_{\Omega,1}e_i,e_j)_{\mathcal E}\\
	&= \lambda_0(A_1e_i,e_j)_{\mathcal E}
	-\lambda_0^2(B_1e_i,e_j)_{\mathcal E}\\
	&= (A_1u_i,u_j)_{\mathcal E}
	-\lambda_0(B_1u_i,u_j)_{\mathcal E}\\
	&= \int_\Omega
	\bigl(S_1(DV)\nabla u_i\bigr)\cdot\nabla u_j\,\mathrm dx
	-\lambda_0\int_\Omega
	(\operatorname{div}V)u_i u_j\,\mathrm dx,
\end{align*}
where the last equality holds by definition of \(A_1\) and \(B_1\). To express the volume integral as a boundary integral, we set
\[
X_{ij}=\bigl(\nabla u_i\cdot\nabla u_j-\lambda_0u_i u_j\bigr)V
-(V\cdot\nabla u_i)\nabla u_j
-(V\cdot\nabla u_j)\nabla u_i.
\]
By direct calculation, we obtain
\[
\operatorname{div}X_{ij}
=\bigl(S_1(DV)\nabla u_i\bigr)\cdot\nabla u_j
-\lambda_0(\operatorname{div}V)u_i u_j,
\]
where $S_1(DV)=(\operatorname{div}V)I_d-DV-(DV)^{\top}$.
On $\partial\Omega$, the Dirichlet condition gives
$\nabla u_k=(\partial_\nu u_k)\nu$ for $k=i,j$. Hence
\[
X_{ij}\cdot\nu=-(V\cdot\nu)\,\partial_\nu u_i\,\partial_\nu u_j
\qquad\text{on }\partial\Omega.
\]
The divergence theorem therefore gives
\[
(H_{\Omega,1}(V))_{ij}
=\int_\Omega\operatorname{div}X_{ij}\,\mathrm dx
=-\int_{\partial\Omega}(V\cdot\nu)\,
\partial_\nu u_i\,\partial_\nu u_j\,\mathrm dS.
\]
Thus $H_{\Omega,1}(V)$ is the Hadamard matrix
in~\eqref{eq:hadamard}, with $h=V\cdot\nu$.

\medskip
\noindent\textbf{Step 4. Estimate for the remainder.} From Step 3, we have
\[
L(tV)=\mu_0\operatorname{Id}_{E_0}+\sum_{k\ge1}t^kL_k(V).
\]
We claim that $L_k(V)$ is a finite linear combination of products
of the following form, possibly integrated over one or more copies of $\Gamma$:
\[
P_0 Y_0 X_1 Y_1\cdots X_qY_qP_0,\qquad
X_r\in\{A_{j_r},B_{j_r}\},\qquad
j_r\ge1,\quad j_1+\cdots+j_q=k.
\]
Each $Y_r$ is a product of
fixed factors from
$\{\operatorname{Id}_{\mathcal E},P_0,B_0,R_0(z):z\in\Gamma\}$. Indeed, since $P(tV)$ commutes with $T_{\Omega}(tV)$ and is a projection, we have
\[
L(tV)=N(tV)^{-1/2}
\left.P_0T_{\Omega}(tV)P(tV)P_0\right|_{E_0}N(tV)^{-1/2}.
\]
We expand the series as in Step~3 and substitute. By direct computation, we obtain the products above for each \(L_k\). As an example, we give the case \(k=2\) explicitly. In Step~3, we get $L_2=P_0T_{\Omega,2}P_0+P_0T_{\Omega,1}P_1P_0$. Substituting the formulas for $T_{\Omega,1}$ and $T_{\Omega,2}$ from Step~3, we obtain
\begin{align*}
	P_0T_{\Omega,2}P_0
	={}& P_0B_2P_0
	-\tfrac12 P_0A_2B_0P_0
	-\tfrac12 P_0B_0A_2P_0
	+\tfrac38 P_0A_1^2B_0P_0
	+\tfrac38 P_0B_0A_1^2P_0 \\
	&-\tfrac12 P_0A_1B_1P_0
	-\tfrac12 P_0B_1A_1P_0
	+\tfrac14 P_0A_1B_0A_1P_0.
\end{align*}
For the second term in $L_2$, the contour formula for $P_1$
gives
\[
P_0T_{\Omega,1}P_1P_0
=\frac{1}{2\pi\mathrm{i}}\int_{\Gamma}
P_0T_{\Omega,1}R_0(z)T_{\Omega,1}R_0(z)P_0\,\mathrm{d}z.
\]
Substituting
$T_{\Omega,1}=B_1
-\tfrac12(A_1B_0+B_0A_1)$
gives the desired result.

We now return to the general case. Each fixed factor $Y_r$ commutes with $P_0$ and acts on $E_0$ as a scalar.
Thus $P_0Y_r=Y_rP_0=c_rP_0$, where $|c_r|\le\|Y_r\|$.
For each $Y_r$, its norm is uniformly bounded
with respect to the contour variables. If $q=1$, then $j_1=k$. We write
$Y_r|_{E_0}=c_r\operatorname{Id}_{E_0}$ for $r=0,1$.
Since $Y_0$ commutes with $P_0$, we have $P_0Y_0=c_0P_0$ and $Y_1P_0=c_1P_0$.
Hence, for $k\ge2$, estimate \eqref{eq:twoends} gives
\[
\|P_0Y_0X_1Y_1P_0\|
=|c_0c_1|\,\|P_0X_1P_0\|
\le |c_0c_1|\,C^k\delta^{k-2}\sigma^2.
\]
If $q\ge2$, we can get, in a similar way,
\[
P_0Y_0X_1Y_1\cdots X_qY_qP_0
=c_0c_q(P_0X_1)Y_1X_2\cdots Y_{q-1}(X_qP_0).
\]
Since $P_0$ and $X_1$ are self-adjoint, $\|P_0X_1\|=\|X_1P_0\|$.
Thus \eqref{eq:oneend} applies to both $P_0X_1$ and $X_qP_0$,
while \eqref{eq:operatorbound} applies to each $X_r$ with $2\le r\le q-1$.
By submultiplicativity of the operator norm,
\begin{align*}
	\|P_0Y_0X_1Y_1\cdots X_qY_qP_0\|
	&\le |c_0c_q|
	\left(\prod_{r=1}^{q-1}\|Y_r\|\right)
	\|P_0X_1\|
	\left(\prod_{r=2}^{q-1}\|X_r\|\right)
	\|X_qP_0\|\\
	&\le |c_0c_q|
	\left(\prod_{r=1}^{q-1}\|Y_r\|\right)
	\bigl(C^{j_1}\delta^{j_1-1}\sigma\bigr)
	\left(\prod_{r=2}^{q-1}C^{j_r}\delta^{j_r}\right)
	\bigl(C^{j_q}\delta^{j_q-1}\sigma\bigr)\\
	&= |c_0c_q|
	\left(\prod_{r=1}^{q-1}\|Y_r\|\right)
	C^k\delta^{k-2}\sigma^2,
\end{align*}
where the last equality uses $j_1+\cdots+j_q=k$.

Denote the coefficient of $t^k$ in $H_\Omega(tV)$ by $H_{\Omega,k}(V,\ldots,V)$, the diagonal value of the symmetric $k$-linear Taylor coefficient. Expanding the inverse $L(tV)^{-1}$ as in Appendix~\ref{app:coefficientweights} expresses this coefficient as a finite sum of products of the $L_j$ of total degree $k$. The preceding estimates then give directly
\[
\|H_{\Omega,k}(V,\ldots,V)\|\le a_k\,\delta^{k-2}\sigma^2 .
\]
As shown in Appendix~\ref{app:coefficientweights}, there exists a
constant $D\ge1$, independent of $k$ and $V$, such that
$a_k\le D^k$ for all $k\ge2$. Thus we have
\[
\|H_{\Omega,k}(V,\ldots,V)\|\le D^k\delta^{k-2}\sigma^2,
\qquad k\ge2.
\]
In particular, $\|H_{\Omega,2}(V,V)\|\le D^2\sigma^2$. Finally, taking $D\delta\le\tfrac12$ and setting $t=1$, we obtain
\[
\|\mathcal R_\Omega(V)\|
\le \sigma^2\sum_{k=3}^{\infty}D^k\delta^{k-2}
= \frac{D^3\delta\sigma^2}{1-D\delta}
\le 2D^3\delta\sigma^2.
\]
This proves~\eqref{eq:clusterexpansion}--\eqref{eq:clusterbounds}.
\end{proof}

\begin{remark}
	Formula~(2.7) gives a second-order expansion of the matrix representing
	the perturbed eigenvalues along each fixed deformation direction.
	However, it does not provide uniform estimates for the second-order
	term and the remainder with respect to the deformation. Such estimates
	are essential for establishing local minimality in a common
	neighborhood rather than along individual deformation directions.
	This is the main purpose of Lemma~\ref{lem:cluster}.
\end{remark}

Lemma~\ref{lem:cluster} gives a matrix decomposition and uniform
estimates in terms of the interior deformation field $V$.
To prove Theorem~\ref{thm:uniform}, we need to express these estimates
in terms of the boundary deformation function $h$ and incorporate
the volume constraint. Moreover, we need to relate the second-order
matrix to the trace formula established in Theorem~\ref{thm:kernel}.
These are the purposes of the following proposition. We would like to point out that, although Lemma~\ref{lem:cluster} is established for general deformation
fields $V$, it is sufficient for our purposes to consider radial
displacements of the boundary of the ball. Indeed, the boundary of
any sufficiently small $C^{2,\alpha}$ deformation of the ball can be
represented as a radial graph over $\mathbb{S}^{d-1}$ after
reparametrization. Since the Dirichlet eigenvalues depend only on the
resulting domain and not on the particular deformation map, this
restriction entails no loss of generality.

\begin{proposition}
Let $h\in C^{2,\alpha}(\mathbb{S}^{d-1})$ be real, have zero mean, and satisfy $\|h\|_{C^{2,\alpha}}<\eta$ for a sufficiently small fixed $\eta\in(0,1)$. The cluster $\lambda_2(\Omega_h),\ldots,\lambda_{d+1}(\Omega_h)$ is represented by a symmetric matrix
\begin{equation}\label{eq:boundarymatrix}
H\left(h\right)=\lambda_*I_d+M\left(h_2\right)+K\left(h,h\right)+R\left(h\right),
\end{equation}
where
\begin{equation}\label{eq:boundaryestimates}
\left\lVert K\left(h,h\right)\right\rVert\le C\left\lVert h\right\rVert_{1/2,*}^2,
\qquad
\left\lVert R\left(h\right)\right\rVert\le C\|h\|_{C^{2,\alpha}}\left\lVert h\right\rVert_{1/2,*}^2.
\end{equation}
The map $K$ extends to a continuous symmetric bilinear map on the zero-mean subspace of $H^{1/2}(\mathbb{S}^{d-1})$. Moreover, we define $Q(f,g)=d^{-1}\operatorname{tr}K(f,g)$ for zero-mean $f,g\in H^{1/2}(\mathbb S^{d-1})$. If $v_0=v_1=v_2=0$, then
\begin{equation}\label{eq:Qtrace}
Q(v,v)
=\frac{2\lambda_*}{s_d}\sum_{n\ge3}A_{d,n}\left\lVert v_n\right\rVert_2^2.
\end{equation}
\end{proposition}
\begin{proof}
We first consider an arbitrary sufficiently small boundary displacement $\psi$, without imposing the volume constraint. Let $E(\psi\theta)$ denote the harmonic extension to $B$ of the vector field $\theta\mapsto\psi(\theta)\theta$. By the trace theorem and the energy-minimizing property of harmonic extension~\cite{McLean}, we have
\begin{equation}\label{eq:harmonicenergy}
\left\lVert DE\left(\psi\theta\right)\right\rVert_{L^2\left(B\right)}
\le C\left\lVert \psi\right\rVert_{1/2,*}.
\end{equation}
Boundary Schauder estimates~\cite{GT} also imply
\begin{equation}\label{eq:harmonicstrong}
\left\lVert DE\left(\psi\theta\right)\right\rVert_{L^\infty\left(B\right)}
\le C\left\lVert \psi\right\rVert_{C^{1,\alpha}}
\le C\left\lVert \psi\right\rVert_{C^{2,\alpha}}.
\end{equation}
For small $\psi$, the map $F_\psi(x)=x+E(\psi\theta)(x)$ is a diffeomorphism of the ball onto the radial domain of radius $1+\psi$. Lemma~\ref{lem:cluster}, together with~\eqref{eq:harmonicenergy}--\eqref{eq:harmonicstrong}, yields a matrix
\[
\widetilde H\left(\psi\right)
=\lambda_*I_d+M\left(\psi\right)+\widetilde K\left(\psi,\psi\right)
+\widetilde R\left(\psi\right),
\]
where
\begin{equation}\label{bound-for-K-R}
\left\lVert \widetilde K\left(\psi,\psi\right)\right\rVert
\le C\left\lVert \psi\right\rVert_{1/2,*}^2,
\qquad
\left\lVert \widetilde R\left(\psi\right)\right\rVert
\le C\left\lVert \psi\right\rVert_{C^{2,\alpha}}\left\lVert \psi\right\rVert_{1/2,*}^2.
\end{equation}
Here $\widetilde H(\psi)=H_B(E(\psi\theta))$. We use the $H_0^1(B)$-orthonormal basis $e_i=u_i/\sqrt{\lambda_*}$ with $u_i$ as in~\eqref{eq:normalization}, so $M(\psi)$ is the Hadamard matrix~\eqref{eq:hadamard} in the corresponding $L^2(B)$-orthonormal basis $(u_i)$.

We now impose the exact volume constraint. Set
\[
\psi(h)=\zeta\left(h\right)\left(1+h\right)-1,
\qquad g_{\mathrm{vol}}\left(h\right)=-\frac{d-1}{2s_d}\left\lVert h\right\rVert_2^2.
\]
Since $1+\psi(h)=\zeta(h)(1+h)$, the radial domain $F_{\psi(h)}(B)$ is exactly $\Omega_h$ in~\eqref{eq:chart}. We therefore take its cluster matrix to be
\[
H(h):=\widetilde H\bigl(\psi(h)\bigr)
     =\lambda_*I_d+M\left(\psi(h)\right)+\widetilde K\left(\psi(h),\psi(h)\right)
     +\widetilde R\left(\psi(h)\right).
\]
Our next task is to express the above equality in the form
of~\eqref{eq:boundarymatrix}. Below we write $\psi=\psi(h)$ and regard $g_{\mathrm{vol}}(h)$ as a constant function on the sphere. Lemma~\ref{lem:volume} gives
\[
\psi=h+g_{\mathrm{vol}}\left(h\right)+e\left(h\right),
\qquad
 e\left(h\right)=\left(\zeta\left(h\right)-1\right)h+r_\zeta\left(h\right).
\]
It is easy to see from~\eqref{eq:hnorm} that $\| c\|_{1/2,*}=s_d^{1/2}\left|c\right|$ for a constant $c$. Thus we have $\|g_{\mathrm{vol}}(h)\|_{1/2,*}\le C\|h\|_{1/2,*}^2$. Moreover, by \eqref{eq:zetaremainder}, we obtain
\begin{equation}\label{bound-e(h)}
\begin{aligned}
	\left\lVert e\left(h\right)\right\rVert_{1/2,*}
	&\le\left|\zeta\left(h\right)-1\right|\|h\|_{1/2,*}+s_d^{1/2}\left|r_\zeta\left(h\right)\right|\\
	&\le C\|h\|_{2}^2\|h\|_{1/2,*}+C\|h\|_{L^\infty}\|h\|_{2}^2\|h\|_{1/2,*}+C\|h\|_{L^\infty}\|h\|_{2}^2\\
	&\le C\|h\|_{C^{2,\alpha}}\|h\|_{1/2,*}^2.
\end{aligned}
\end{equation}
In the last step, we used $\|h\|_{2}\le\|h\|_{1/2,*}\le C\|h\|_{C^{2,\alpha}}\le C$ and $\|h\|_{L^\infty}\le\|h\|_{C^{2,\alpha}}$. Together with the preceding estimate for $g_{\mathrm{vol}}(h)$, this gives
\begin{equation}\label{boubd-for-psi-h}
\|\psi-h\|_{1/2,*}\le C\|h\|_{1/2,*}^2+C\|h\|_{C^{2,\alpha}} \|h\|_{1/2,*}^2\le C\|h\|_{1/2,*}^2.
\end{equation}
The triangle inequality also gives
\begin{equation}\label{bound-for-psi}
\begin{aligned}
\|\psi\|_{1/2,*}
&\le\|h\|_{1/2,*}+\|\psi-h\|_{1/2,*}\\
&\le\|h\|_{1/2,*}+C\|h\|_{1/2,*}^2\le C\|h\|_{1/2,*},
\end{aligned}
\end{equation}
where we used $\|h\|_{C^{2,\alpha}}\le \eta<1$. Using $\psi=\zeta(h)h+\zeta(h)-1$ and \eqref{eq:zetaremainder}, we obtain
\begin{equation}\label{eq:psibounds}
\left\lVert \psi\right\rVert_{C^{2,\alpha}}\le C\|h\|_{C^{2,\alpha}}.
\end{equation}

We first consider $\widetilde K\left(\psi(h),\psi(h)\right)$. We may substitute \(\psi(h)=h+g_{\mathrm{vol}}+e\) and expand directly into six terms, but this is rather tedious. After some analysis, we instead examine the expression with the second-order $\widetilde K\left(h,h\right)$ term subtracted. From \eqref{bound-for-K-R}, polarization of the quadratic bound for $\widetilde K$ yields $\| \widetilde K(f,g) \|
\le C\left\lVert f\right\rVert_{1/2,*}\left\lVert g\right\rVert_{1/2,*}$. Consequently, by \eqref{boubd-for-psi-h} and \eqref{bound-for-psi},
\begin{equation}\label{K-K}
\begin{aligned}
	\left\lVert \widetilde K\left(\psi,\psi\right)-\widetilde K\left(h,h\right)\right\rVert
	&=\left\lVert \widetilde K\left(\psi-h,\psi+h\right)\right\rVert\\
	&\le C\left\lVert \psi-h\right\rVert_{1/2,*}
	\left(\left\lVert \psi\right\rVert_{1/2,*}+\left\lVert h\right\rVert_{1/2,*}\right)\\
	&\le C\|h\|_{1/2,*}^3\le C\|h\|_{C^{2,\alpha}}\|h\|_{1/2,*}^2.
\end{aligned}
\end{equation}
The above relation shows that the order of $\widetilde K\left(\psi,\psi\right)-\widetilde K\left(h,h\right)$ is at least 3. For the linear term $M\left(\psi(h)\right)$, we have $M\left(\psi(h)\right)=M\left(h\right)+M\left(g_{\mathrm{vol}}(h)\right)+M\left(e(h)\right)$. For zero-mean $h$, $M(h)=M(h_2)$ by Lemma~\ref{lem:firstmatrix}. By \eqref{eq:hadamard-norm}, \eqref{bound-e(h)} and \eqref{constant-M(c)}, we obtain
\begin{equation}\label{Meh}
\left\lVert M(e(h))\right\rVert\le C\|h\|_{C^{2,\alpha}}\|h\|_{1/2,*}^2,
\end{equation}
\begin{equation}\label{gvol}
M\bigl(g_{\mathrm{vol}}(h)\bigr)
=-2\lambda_*\left(-\frac{d-1}{2s_d}\|h\|_2^2\right)I_d
=\frac{\lambda_*(d-1)}{s_d}\|h\|_2^2 I_d.
\end{equation}
By~\eqref{gvol}, $M\bigl(g_{\mathrm{vol}}(h)\bigr)$ is quadratic
in $h$, so we can collect it with the other second-order terms:
\[
K(h,h)=\widetilde K(h,h)+M\bigl(g_{\mathrm{vol}}(h)\bigr).
\]
By~\eqref{K-K} and \eqref{Meh}, the terms of order at least $3$ form the remainder
\[
R(h)=M(e(h))+\widetilde K(\psi,\psi)-\widetilde K(h,h)+\widetilde R(\psi).
\]
Thus we obtain
\begin{equation}
	\begin{aligned}
		H(h)
		&= \lambda_* I_d + M\big(\psi(h)\big) + \widetilde K\big(\psi(h),\psi(h)\big) + \widetilde R\big(\psi(h)\big) \\
		&= \lambda_* I_d + M(h) + M\big(g_{\mathrm{vol}}(h)\big) + M\big(e(h)\big) + \widetilde K(h,h) + \widetilde K(\psi,\psi) - \widetilde K(h,h) + \widetilde R(\psi) \\
		&= \lambda_* I_d + M(h_2) + K(h,h) + R(h).
	\end{aligned}
\end{equation}
This proves~\eqref{eq:boundarymatrix}. The preceding estimates \eqref{bound-for-K-R}, \eqref{eq:psibounds}--\eqref{gvol} give \eqref{eq:boundaryestimates}.

The bound $\|K(f,g)\|\le C\|f\|_{1/2,*}\|g\|_{1/2,*}$ and the density of smooth zero-mean functions imply that $K$ extends uniquely to a continuous symmetric bilinear map on the
zero-mean subspace of $H^{1/2}(\mathbb S^{d-1})$. Let $v^{(N)}=\sum_{n=3}^{N}v_n$ be a finite sum of spherical
harmonics. Since $M(v^{(N)})=0$, the expansion of
$H(\varepsilon v^{(N)})$ is
\[
\begin{aligned}
	\frac{1}{d}\sum_{j=2}^{d+1}\lambda_j(\Omega_{\varepsilon v^{(N)}})
	=\frac{1}{d}\operatorname{tr}H(\varepsilon v^{(N)})
	=\lambda_*+\varepsilon^2 Q(v^{(N)},v^{(N)})
	+O_{v^{(N)}}(|\varepsilon|^3).
\end{aligned}
\]
For the same $v^{(N)}$,
Theorem~\ref{thm:kernel} gives
\[
\frac{1}{d}\sum_{j=2}^{d+1}\lambda_j(\Omega_{\varepsilon v^{(N)}})
=\lambda_*+\varepsilon^2\frac{2\lambda_*}{s_d}
\sum_{n=3}^{N}A_{d,n}\|v_n\|_2^2
+o_{v^{(N)}}(\varepsilon^2).
\]
Comparing the coefficients of $\varepsilon^2$, we have
\[
Q(v^{(N)},v^{(N)})=\frac{2\lambda_*}{s_d}
\sum_{n=3}^{N}A_{d,n}\|v_n\|_2^2.
\]
Now let $v\in H^{1/2}(\mathbb S^{d-1})$ satisfy $v_0=v_1=v_2=0$. By~\eqref{eq:hnorm},
\[
\|v-v^{(N)}\|_{1/2,*}^2
=\sum_{n>N}(1+n)\|v_n\|_2^2\longrightarrow 0.
\]
Thus, the bilinear bound for $Q$ gives
\[
\bigl|Q(v,v)-Q(v^{(N)},v^{(N)})\bigr|
\le C\|v-v^{(N)}\|_{1/2,*}
\bigl(\|v\|_{1/2,*}+\|v^{(N)}\|_{1/2,*}\bigr)
\longrightarrow 0.
\]
We take the limit in the finite sum identity for \(v^{(N)}\) and prove \eqref{eq:Qtrace}.
\end{proof}

\begin{proof}[Proof of Theorem~\ref{thm:uniform}]
Since $h_0=h_1=0$, we have $h=h_2+h_{\ge3}$. By orthogonality in~\eqref{eq:hnorm},
\[
\|h\|_{1/2,*}^2=\|h_2\|_{1/2,*}^2+\|h_{\ge3}\|_{1/2,*}^2.
\]
We denote the eigenvalue difference by $D_h=\lambda_{d+1}\left(\Omega_h\right)-\lambda_*$. Here $\lambda_{d+1}(\Omega_h)$ is the largest eigenvalue of $H(h)$ when $\|h\|_{C^{2,\alpha}}$ is sufficiently small. We first use the trace of $H(h)$ to obtain a lower bound for $D_h$. Equations~\eqref{eq:coercivity} and~\eqref{eq:Qtrace} imply
\[
Q(h_{\ge3},h_{\ge3})\ge\kappa_d\|h_{\ge3}\|_{1/2,*}^2,
\qquad\kappa_d:=\frac{2\lambda_*c_d}{s_d}>0.
\]
By $|Q(f,g)|\le C_0\|f\|_{1/2,*}\|g\|_{1/2,*}$, we obtain
\begin{align*}
Q(h,h)
&=Q(h_{\ge3},h_{\ge3})+2Q\left(h_2,h_{\ge3}\right)+Q(h_2,h_2)\\
&\ge\kappa_d\|h_{\ge3}\|_{1/2,*}^2-2C_0\|h_2\|_{1/2,*}\|h_{\ge3}\|_{1/2,*}-C_0\|h_2\|_{1/2,*}^2.
\end{align*}
Using Young's inequality in the explicit form
\[
2C_0\|h_2\|_{1/2,*}\|h_{\ge3}\|_{1/2,*}\le\frac{\kappa_d}{2}\|h_{\ge3}\|_{1/2,*}^2+
\frac{2C_0^2}{\kappa_d}\|h_2\|_{1/2,*}^2,
\]
then we have
\[
Q(h,h)\ge\frac{\kappa_d}{2}\|h_{\ge3}\|_{1/2,*}^2-C_1\|h_2\|_{1/2,*}^2,
\qquad C_1=C_0+\frac{2C_0^2}{\kappa_d}.
\]
Since $\operatorname{tr} M(h_2)=0$ and $\left|d^{-1}\operatorname{tr} R(h)\right|\le\left\lVert R(h)\right\rVert$, we obtain
\[
\begin{aligned}
	D_h
	&=\lambda_{\max}\bigl(M(h_2)+K(h,h)+R(h)\bigr)\\
	&\ge \frac{1}{d}\operatorname{tr}\bigl(K(h,h)+R(h)\bigr)\\
	&\ge \frac{\kappa_d}{2}\|h_{\ge3}\|_{1/2,*}^2
	-C_1\|h_2\|_{1/2,*}^2
	-C_R\|h\|_{C^{2,\alpha}}
	\left(\|h_2\|_{1/2,*}^2+\|h_{\ge3}\|_{1/2,*}^2\right).
\end{aligned}
\]
Choose $\eta\in(0,1)$ small enough that $C_R\eta\le\kappa_d/4$. Then
\begin{equation}\label{eq:quadraticbound}
D_h\ge\frac{\kappa_d}{4}\|h_{\ge3}\|_{1/2,*}^2-C_2\|h_2\|_{1/2,*}^2,
\qquad C_2=C_1+C_R.
\end{equation}

We next use the largest eigenvalue of the first-order matrix. The variational characterization of the largest matrix eigenvalue implies $\lambda_{\max}\left(S+E\right)\ge\lambda_{\max}\left(S\right)-\left\lVert E\right\rVert$ for symmetric $S,E$. Applying this to~\eqref{eq:boundarymatrix}, and using~\eqref{eq:positivefirst} and~\eqref{eq:boundaryestimates}, we get
\begin{equation}\label{eq:linearbound}
D_h\ge m_d\|h_2\|_{1/2,*}-C_3\left(\|h_2\|_{1/2,*}^2+\|h_{\ge3}\|_{1/2,*}^2\right),
\end{equation}
where $C_3>0$ is independent of $h$. We choose $0<\tau\le\frac{\kappa_d}{8C_3+2\kappa_d}$.
Then $(1-\tau)\frac{\kappa_d}{4}-\tau C_3\ge\frac{\kappa_d}{8}$. Taking a convex combination of~\eqref{eq:quadraticbound} and~\eqref{eq:linearbound} yields
\begin{align*}
	D_h
	&\ge (1-\tau)\left(
	\frac{\kappa_d}{4}\|h_{\ge3}\|_{1/2,*}^2
	-C_2\|h_2\|_{1/2,*}^2\right)+\tau\left[
	m_d\|h_2\|_{1/2,*}
	-C_3\left(\|h_2\|_{1/2,*}^2
	+\|h_{\ge3}\|_{1/2,*}^2\right)\right]\\
	&=\tau m_d\|h_2\|_{1/2,*}
	+\left((1-\tau)\frac{\kappa_d}{4}-\tau C_3\right)
	\|h_{\ge3}\|_{1/2,*}^2
	-\bigl((1-\tau)C_2+\tau C_3\bigr)\|h_2\|_{1/2,*}^2\\
	&\ge \tau m_d\|h_2\|_{1/2,*}
	+\frac{\kappa_d}{8}\|h_{\ge3}\|_{1/2,*}^2
	-C_4\|h_2\|_{1/2,*}^2,
\end{align*}
where $C_4=(1-\tau)C_2+\tau C_3$. Using the continuous embedding $\|h_2\|_{1/2,*}\le\|h\|_{1/2,*}\le C_5\|h\|_{C^{2,\alpha}}$, we shrink the common neighborhood until $C_4C_5\|h\|_{C^{2,\alpha}}\le\frac{\tau m_d}{2}$. Then we can absorb $C_4\|h_2\|_{1/2,*}^2$ into half of the linear term and obtain
\begin{align*}
D_h&\ge\frac{\tau m_d}{2}\|h_2\|_{1/2,*}+\frac{\kappa_d}{8}\|h_{\ge3}\|_{1/2,*}^2\ge c\left(\|h_2\|_{1/2,*}+\|h_{\ge3}\|_{1/2,*}^2\right).
\end{align*}
This proves the lower bound in~\eqref{eq:mainuniform}.

For the upper bound, boundedness of the degree-two linear map and~\eqref{eq:boundaryestimates} imply
\begin{align*}
D_h&\le\left\lVert M\left(h_2\right)\right\rVert+\left\lVert K\left(h,h\right)\right\rVert+\left\lVert R\left(h\right)\right\rVert\\
&\le C_6\|h_2\|_{1/2,*}+C_7\left(\|h_2\|_{1/2,*}^2+\|h_{\ge3}\|_{1/2,*}^2\right).
\end{align*}
Taking $\eta$ smaller if necessary, we have
$\|h_2\|_{1/2,*}\le C_5\eta\le 1$ and $\|h_2\|_{1/2,*}^2\le \|h_2\|_{1/2,*}$. Hence, we obtain the upper bound in~\eqref{eq:mainuniform}. Finally,
\[
\|h_2\|_{1/2,*}+\|h_{\ge3}\|_{1/2,*}^2\ge \|h_2\|_{1/2,*}^2+\|h_{\ge3}\|_{1/2,*}^2=\|h\|_{1/2,*}^2
\]
proves~\eqref{eq:quadraticlocal}. Every choice of a constant or neighborhood size depends only on the fixed dimension, H\"older exponent, and ball cluster.
\end{proof}

\paragraph{A geometric form of the local estimate.}
We translate the boundary norm into a measure of the difference between sets. For a measurable set $\Omega$ of finite positive volume, its Fraenkel asymmetry is
\begin{equation}\label{eq:asymmetry}
\mathcal A\left(\Omega\right)
=\inf_{y\in\mathbb{R}^d}\frac{\left|\Omega\mathbin\triangle\left(y+B_\Omega\right)\right|}{\left|\Omega\right|},
\end{equation}
where $B_\Omega=B_{(|\Omega|/\omega_d)^{1/d}}$ is the ball centered at the origin having volume $\left|\Omega\right|$ and $\triangle$ denotes symmetric difference. This quantity is invariant under translations and dilations. It is the geometric distance used in the quantitative Faber--Krahn inequality~\cite{BDPV}.

\begin{corollary}\label{cor:asymmetry}
Under the assumptions of Theorem~\ref{thm:uniform}, we have
\[
\lambda_{d+1}\left(\Omega_h\right)-\lambda_*
\ge c\mathcal A\left(\Omega_h\right)^2.
\]
\end{corollary}
\begin{proof}
We set $\rho=\zeta(h)(1+h)$. It is easy to see that
\[
\begin{aligned}
	|\Omega_h\triangle B|
	&=\int_{\mathbb S^{d-1}}
	\int_{\min\{1,\rho(\theta)\}}^{\max\{1,\rho(\theta)\}}
	r^{d-1}\,\mathrm{d}r\,\mathrm{d}S(\theta)=\frac{1}{d}\int_{\mathbb S^{d-1}}
	\bigl|\rho(\theta)^d-1\bigr|\,\mathrm{d}S(\theta).
\end{aligned}
\]
For small $\left\lVert h\right\rVert_{L^\infty}$, both $\rho$ and $\zeta(h)$ stay in a fixed compact subinterval of $(0,\infty)$. By Lemma~\ref{lem:volume}, we have
\begin{align*}
\left|\Omega_h\triangle B\right|
&\le C\int_{\mathbb{S}^{d-1}}\left|\rho-1\right|\,\mathrm{d} S\\
&\le C\left(\left\lVert h\right\rVert_{L^1\left(\mathbb{S}^{d-1}\right)}+\left\lVert h\right\rVert_2^2\right)\\
&\le C\left\lVert h\right\rVert_2\le C\left\lVert h\right\rVert_{1/2,*}.
\end{align*}
Combining with~\eqref{eq:quadraticlocal} and~\eqref{eq:asymmetry}, we prove the claim.
\end{proof}

\begin{remark}
Corollary~\ref{cor:asymmetry} has a geometric right-hand side, but its hypotheses still require a small $C^{2,\alpha}$ perturbation. We do not assert a neighborhood theorem for arbitrary sets of small symmetric difference.
\end{remark}

\paragraph{The two-ball configuration for $\lambda_{d+2}$.}
We finally apply the single ball result to Open Problem~11.5. We first give the radii of the reference balls. Given a total volume $\mathcal V>0$, we aim to construct two balls $B_{R_1}$ and $B_{R_2}$ with radii $R_1>R_2$ such that the \((d+1)\)-th eigenvalue of the larger ball coincides with the first eigenvalue of the smaller one. Thus let
\[
\lambda_{d+1}\left(B_{R_1}\right)=\frac{\lambda_*}{R_1^2}=\lambda_1\left(B_{R_2}\right)=\frac{\lambda_1\left(B\right)}{R_2^2}:=\Lambda_*.
\]
Then, the volume constraint $\omega_d(R_1^d+R_2^d)=\mathcal V$ gives
$$
\Lambda_*=\left(\frac{\omega_d\left(\lambda_*^{d/2}+\lambda_1\left(B\right)^{d/2}\right)}{\mathcal V}\right)^{2/d}.
$$
Thus, $B_{R_1}$ has exactly one eigenvalue below $\Lambda_*$, namely $\lambda_1(B)/R_1^2$. Its next $d$ eigenvalues equal $\Lambda_*$. The first eigenvalue of $B_{R_2}$ is $\Lambda_*$, and its second eigenvalue is $\lambda_*/R_2^2>\Lambda_*$. After translating the balls so that they are disjoint, we retain the notation $B_{R_1}$ and $B_{R_2}$. For their union, we have
\[
\lambda_1(B_{R_1}\cup B_{R_2})
<\lambda_2(B_{R_1}\cup B_{R_2})
=\cdots
=\lambda_{d+2}(B_{R_1}\cup B_{R_2})
=\Lambda_*
<\lambda_{d+3}(B_{R_1}\cup B_{R_2}).
\]
The next corollary gives a lower bound for the eigenvalue increase in terms of changes in the shapes and volumes of the two components. The normalization of the first component uses the chart fixed in~\eqref{eq:chart}.

\begin{corollary}\label{cor:twoball}
Fix $d\ge2$, $0<\alpha<1$, and $\mathcal V>0$. Let $\Omega_1,\Omega_2$ be disjoint smooth domains having total volume $\mathcal V$. Suppose that $\left|\Omega_1\right|$ is close to $|B_{R_1}|$ and that, after  rescalings to volume $\omega_d$, $\Omega_1,\Omega_2$ are small $C^{2,\alpha}$ perturbations of $B$. We write the normalized first component as $\Omega_h$ in the chart $h_0=h_1=0$. Here $|\Omega_h|=\omega_d$. Then, in one fixed neighborhood of this configuration,
\begin{equation}\label{eq:twoballstability}
\lambda_{d+2}\left(\Omega_1\cup \Omega_2\right)-\Lambda_*
\ge c\left(\Bigl\lvert \lvert\Omega_1\rvert - \lvert B_{R_1}\rvert \Bigr\rvert+\mathscr D\left(h\right)+\mathcal A\left(\Omega_2\right)^2\right).
\end{equation}
The constant $c>0$ depends only on $d$, $\alpha$ and $\mathcal V$.
\end{corollary}
\begin{proof}
Since the eigenvalues depend continuously on smooth deformations, the $(d+2)$-th eigenvalue of the perturbed domain is
\begin{equation}\label{eq:ordering}
\lambda_{d+2}\left(\Omega_1\cup \Omega_2\right)
=\max\left\{\lambda_{d+1}\left(\Omega_1\right),\lambda_1\left(\Omega_2\right)\right\}.
\end{equation}
Let $C_1=\lambda_*\omega_d^{2/d}$ and $C_2=\lambda_1(B)\omega_d^{2/d}$. By Theorem~\ref{thm:uniform} and \eqref{scaling-laws}, we obtain
\begin{align*}
	\lambda_{d+1}(\Omega_1)
	=\left(\frac{\omega_d}{|\Omega_1|}\right)^{2/d}
	\lambda_{d+1}(\Omega_h)
	\ge C_1\bigl(1+c_1\mathscr D(h)\bigr)|\Omega_1|^{-2/d}
\end{align*}
for some constant $c_1>0$. For $\lambda_1\left(\Omega_2\right)$, we use the quantitative Faber--Krahn inequality~\cite{BDPV} 
\[
\lambda_1(\Omega)|\Omega|^{2/d}
-\lambda_1(B)\omega_d^{2/d}
\ge c_d^{\mathrm{FK}}\mathcal A(\Omega)^2,
\qquad c_d^{\mathrm{FK}}>0.
\]
Since $|\Omega_2|=\mathcal V-|\Omega_1|$, it follows that
\[
\lambda_1(\Omega_2)
\ge C_2\bigl(1+c_2\mathcal A(\Omega_2)^2\bigr)
\bigl(\mathcal V-|\Omega_1|\bigr)^{-2/d}
\]
for some constant $c_2>0$.

Let $U=C_1\bigl(1+c_1\mathscr D(h)\bigr)$ and $W=C_2\bigl(1+c_2\mathcal A(\Omega_2)^2\bigr)$. For \(U,W>0\), the function \(s\mapsto Us^{-2/d}\) is strictly decreasing on \((0,\mathcal V)\), whereas \(s\mapsto W(\mathcal V-s)^{-2/d}\) is strictly increasing. Their maximum is minimized at their unique intersection,
where
\[
s=\frac{\mathcal V U^{d/2}}{U^{d/2}+W^{d/2}}.
\]
Substituting $s=\frac{\mathcal V U^{d/2}}{U^{d/2}+W^{d/2}}$, we get
\begin{equation}\label{eq:allocation}
	\min_{0<s<\mathcal V}
	\max\bigl\{Us^{-2/d},W(\mathcal V-s)^{-2/d}\bigr\}
	=\Psi(U,W),
\end{equation}
where
\[
\Psi(U,W)
=\mathcal V^{-2/d}\bigl(U^{d/2}+W^{d/2}\bigr)^{2/d}.
\]
Thus, we obtain from~\eqref{eq:ordering} and~\eqref{eq:allocation} that
\[
\lambda_{d+2}(\Omega_1\cup\Omega_2)
\ge\max\bigl\{U|\Omega_1|^{-2/d},
W(\mathcal V-|\Omega_1|)^{-2/d}\bigr\}
\ge\Psi(U,W).
\]
When $h=0$, $\Psi(C_1,C_2)=\Lambda_*$. By direct calculation, we have
\begin{align*}
	\partial_U\Psi
	=\mathcal V^{-2/d}
	\bigl(U^{d/2}+W^{d/2}\bigr)^{2/d-1}U^{d/2-1},\qquad
	\partial_W\Psi
	=\mathcal V^{-2/d}
	\bigl(U^{d/2}+W^{d/2}\bigr)^{2/d-1}W^{d/2-1}.
\end{align*}
It is easy to see that $\partial_U\Psi$ and $\partial_W\Psi$ have a positive lower bound $c_3$ on a fixed compact
rectangle in $(0,\infty)^2$ containing $(C_1,C_2)$ and all admissible
pairs $(U,W)$. Since $U\ge C_1$ and $W\ge C_2$, integration gives
\begin{equation}\label{shapes-change}
\begin{aligned}
	\lambda_{d+2}(\Omega_1\cup\Omega_2)-\Lambda_*
	&\ge \Psi(U,W)-\Psi(C_1,C_2)\\
	&=\int_{C_1}^{U}\partial_U\Psi(s,C_2)\,\mathrm{d}s
	+\int_{C_2}^{W}\partial_W\Psi(U,t)\,\mathrm{d}t\\
	&\ge c_3\bigl((U-C_1)+(W-C_2)\bigr)\\
	&\ge c_4\bigl(\mathscr D(h)+\mathcal A(\Omega_2)^2\bigr).
\end{aligned}
\end{equation}

It remains to estimate the change in the volume allocation. At the
reference configuration,
\[
C_1|B_{R_1}|^{-2/d}
=C_2\bigl(\mathcal V-|B_{R_1}|\bigr)^{-2/d}
=\Lambda_*.
\]
If $|\Omega_1|\le |B_{R_1}|$, we use $U\ge C_1$ to obtain
\begin{align*}
	\lambda_{d+2}(\Omega_1\cup\Omega_2)-\Lambda_*
	&\ge C_1|\Omega_1|^{-2/d}-C_1|B_{R_1}|^{-2/d}\\
	&=\frac{2C_1}{d}\int_{|\Omega_1|}^{|B_{R_1}|}
	s^{-2/d-1}\,\mathrm{d}s\\
	&\ge c_5\bigl(|B_{R_1}|-|\Omega_1|\bigr).
\end{align*}
If $|\Omega_1|\ge |B_{R_1}|$, we use $W\ge C_2$ to obtain
\begin{align*}
	\lambda_{d+2}(\Omega_1\cup\Omega_2)-\Lambda_*
	&\ge C_2\bigl(\mathcal V-|\Omega_1|\bigr)^{-2/d}
	-C_2\bigl(\mathcal V-|B_{R_1}|\bigr)^{-2/d}\\
	&=\frac{2C_2}{d}\int_{|B_{R_1}|}^{|\Omega_1|}
	(\mathcal V-s)^{-2/d-1}\,\mathrm{d}s\\
	&\ge c_6\bigl(|\Omega_1|-|B_{R_1}|\bigr).
\end{align*}
The constants are uniform because $|\Omega_1|\in(0,\mathcal V)$. Thus
\begin{equation}\label{change-volume}
\lambda_{d+2}(\Omega_1\cup\Omega_2)-\Lambda_*
\ge c_7\bigl|\,|\Omega_1|-|B_{R_1}|\,\bigr|,
\qquad c_7=\min\{c_5,c_6\}.
\end{equation}
By \eqref{shapes-change} and \eqref{change-volume}, we prove~\eqref{eq:twoballstability}.
\end{proof}

\begin{remark}
	Under the assumptions of Corollary~\ref{cor:twoball},
	$\lambda_{d+2}(\Omega_1\cup\Omega_2)=\Lambda_*$ if and only if
	$\Omega_1$ and $\Omega_2$ are balls of radii $R_1$ and $R_2$, respectively.
\end{remark}

\begin{remark}
	Corollary~\ref{cor:twoball} gives an affirmative answer to the local
	version of Open Problem~11.5. More precisely, the disjoint union of two
	balls $B_{R_1}$ and $B_{R_2}$ with radii satisfying
	\[
	\frac{R_1}{R_2}=\frac{j_{d/2,1}}{j_{d/2-1,1}}
	\]
	is a local minimizer of $\lambda_{d+2}$ under sufficiently small
	$C^{2,\alpha}$ deformations preserving the two-component structure
	and the total volume.
\end{remark}

\appendix

\section{Bessel zero and boundary operator estimates}\label{app:bessel}
We give the estimates used for the reduced boundary operator introduced in Section~\ref{subsec:reduced}.

\begin{lemma}\label{lem:bracket}
For $d\ge2$, we have
\begin{equation}
3\left(d+2\right)<\lambda_*<\left(d+4\right)\left(d+2\right).
\end{equation}
\end{lemma}
\begin{proof}
Set $\mu=d/2\ge1$. The bound of Lorch~\cite{Lorch} gives
\[
j_{\mu,1}^2>\left(\mu+1\right)\left(\mu+5\right)\ge6\left(\mu+1\right)=3\left(d+2\right).
\]
The classical Rayleigh sums~\cite[equation (13)]{HS} satisfy
\[
\sum_{k=1}^{\infty}\frac{1}{j_{\mu,k}^{2}}
=\frac{1}{4\left(\mu+1\right)},
\qquad
\sum_{k=1}^{\infty}\frac{1}{j_{\mu,k}^{4}}
=\frac{1}{16\left(\mu+1\right)^2\left(\mu+2\right)}.
\]
Since
\[
\sum_{k=1}^{\infty}\frac{1}{j_{\mu,k}^{4}}
<\frac{1}{j_{\mu,1}^{2}}
\sum_{k=1}^{\infty}\frac{1}{j_{\mu,k}^{2}}.
\]
we obtain
\[
j_{\mu,1}^2<4\left(\mu+1\right)\left(\mu+2\right)=\left(d+2\right)\left(d+4\right).
\]
\end{proof}

\begin{lemma}\label{lem:tcoeff}
The numbers $t_\ell$, $\ell\ge2$, are strictly increasing. In particular, $t_2=-d$, $t_3=\lambda_*/(d+2)-(d+1)$ and
\begin{equation}\label{eq:lowcoeff}
t_4=\left(d+2\right)\left(\frac{\lambda_*}{\left(d+4\right)\left(d+2\right)-\lambda_*}-1\right).
\end{equation}
For fixed $d$, we also have
\begin{equation}\label{eq:tasymptotic}
t_\ell=\ell+O_d\left(\ell^{-1}\right)\qquad \ell\to\infty.
\end{equation}
Consequently, $T$ is bounded from $\mathcal X\subset H^{1/2}(\mathbb{S}^{d-1})$ to $H^{-1/2}(\mathbb{S}^{d-1})$.
\end{lemma}
\begin{proof}
We first prove the strict monotonicity of $t_\ell$. Let $p=\ell+d/2-1$ and $x=j_{d/2,1}$. By the Mittag--Leffler expansion for Bessel functions, we have
\[
\frac{J_{p+1}(x)}{J_p(x)}
=\sum_{k=1}^{\infty}\frac{2x}{j_{p,k}^2-x^2}.
\]
Combining this identity with the recurrence relation
\[
J_p'(x)=\frac{p}{x}J_p(x)-J_{p+1}(x),
\]
we obtain
\[
t_\ell
=\ell-j_{d/2,1}
\frac{J_{\ell+d/2}(j_{d/2,1})}
{J_{\ell+d/2-1}(j_{d/2,1})}
=\ell-\sum_{k=1}^{\infty}
\frac{2\lambda_*}{j_{\ell+d/2-1,k}^2-\lambda_*}.
\]
Recall that, for each fixed $k\ge1$, the zero $j_{p,k}$
is strictly increasing with respect to $p\ge0$. Thus, for
$m>\ell\ge2$, we have
\[
j_{m+d/2-1,k}>j_{\ell+d/2-1,k}
\ge j_{\ell+d/2-1,1}>j_{d/2,1}.
\]
Consequently,
\[
\begin{aligned}
	t_m-t_\ell
	&=m-\ell+2\lambda_*\sum_{k=1}^{\infty}
	\left(
	\frac{1}{j_{\ell+d/2-1,k}^2-\lambda_*}
	-\frac{1}{j_{m+d/2-1,k}^2-\lambda_*}
	\right)
	>m-\ell>0.
\end{aligned}
\]
This proves the strict monotonicity of $t_\ell$ for $\ell\ge2$. To calculate the first few values, we write $\mu=d/2$ and use
\[
J'_p\left(x\right)=J_{p-1}\left(x\right)-\frac{p}{x}J_p\left(x\right).
\]
Then $t_\ell$ becomes
\[
t_\ell=j_{d/2,1}\frac{J_{\mu+\ell-2}\left(j_{d/2,1}\right)}{J_{\mu+\ell-1}\left(j_{d/2,1}\right)}
-\left(\ell+d-2\right).
\]
For $\ell=2$, the numerator is zero, so $t_2=-d$. For $\ell=3$, we use $J_{\mu+2}(j_{d/2,1})=(d+2)J_{\mu+1}(j_{d/2,1})/j_{d/2,1}$ and obtain $t_3=\frac{\lambda_*}{d+2}-\left(d+1\right)$. A second recurrence yields
\[
J_{\mu+3}\left(j_{d/2,1}\right)
=\frac{\left(d+4\right)\left(d+2\right)-\lambda_*}{\lambda_*}J_{\mu+1}\left(j_{d/2,1}\right),
\]
and hence
\[
t_4
=\frac{\left(d+2\right)\lambda_*}{\left(d+4\right)\left(d+2\right)-\lambda_*}-\left(d+2\right).
\]
This proves~\eqref{eq:lowcoeff}.

We also use the standard fixed-argument large-order estimate from the Bessel series~\cite[Section~10.19]{DLMF}:
\[
J_p\left(j_{d/2,1}\right)
=\frac{(j_{d/2,1}/2)^p}{\Gamma\left(p+1\right)}\left(1+O_d\left(p^{-1}\right)\right).
\]
It follows that
\[
\frac{J_{p+1}\left(j_{d/2,1}\right)}{J_p\left(j_{d/2,1}\right)}
=\frac{j_{d/2,1}}{2\left(p+1\right)}\left(1+O_d\left(p^{-1}\right)\right).
\]
The alternative derivative recurrence gives
\[
t_\ell
=\ell-j_{d/2,1}\frac{J_{\ell+d/2}\left(j_{d/2,1}\right)}{J_{\ell+d/2-1}\left(j_{d/2,1}\right)}
=\ell+O_d\left(\ell^{-1}\right),
\]
which proves~\eqref{eq:tasymptotic}. Thus $\left|t_\ell\right|\le C(1+\ell)$, also for the admissible degree $\ell=0$. 

We denote by $\|\cdot\|_{-1/2,*}$ the dual norm associated
with the norm defined in~\eqref{eq:hnorm}. By the orthogonality
of spherical harmonics, we have
\[
\left\lVert g\right\rVert_{-1/2,*}^2
=\sum_{\ell=0}^{\infty}
\frac{\left\lVert g_\ell\right\rVert_2^2}{1+\ell}.
\]
Since $Tf=\sum_{\ell\ne1}t_\ell f_\ell$ and
$|t_\ell|\le C(1+\ell)$ for $\ell\ne1$, we obtain
\[
\begin{aligned}
	\left\lVert Tf\right\rVert_{-1/2,*}^2
	=\sum_{\ell\ne1}
	\frac{t_\ell^2\left\lVert f_\ell\right\rVert_2^2}{1+\ell}
	\le C^2\left\lVert f\right\rVert_{1/2,*}^2.
\end{aligned}
\]
This proves the boundedness of $T$ from $\mathcal X$
to $H^{-1/2}(\mathbb{S}^{d-1})$.
\end{proof}

\section{An exponential bound for the coefficients}
\label{app:coefficientweights}

In this section, we prove $a_k\le D^k$ used in Step~4. Recall that $a_k$
is the total weight of the terms in the expansion of $H_{\Omega,k}(V,\ldots,V)$. By Step~2, each $A_i$ or $B_i$ of degree
$i\ge1$ has weight at most $C^i$. We first give two basic estimates for the coefficients
of power series.

\medskip
\noindent\textbf{Rule 1: Products of two or three series.}
Consider a product of $r$ series, where $r=2$ or $r=3$. We suppose
that the weight of the constant coefficient is at most
\(\beta_0\ge1\) and the weight of the coefficient of degree \(i\ge1\) is at
most \(\beta_1^i\). To get a term of degree $j$, we choose a term from
each factor such that the sum of their degrees is $j$. For
$i_1+\cdots+i_r=j$, the positive degree coefficients contribute at
most $\beta_1^{i_1+\cdots+i_r}=\beta_1^j$, and the constant coefficients
contribute at most $\beta_0^r$. The number of combinations is $j+1$ for two
factors and $\binom{j+2}{2}$ for three factors. Since $j+1\le 2^j$ and $\binom{j+2}{2}\le 3^j, j\ge1$, the total weight at degree $j$ is at most $\beta_0^r(r\beta_1)^j$. In particular, the total weight is bounded by $C^j$ for some constant $C$.

\medskip
\noindent\textbf{Rule 2: powers of a series.}
Let $F(t)=\sum_{i\ge1}t^iF_i$, where $F_i$ has total weight at
most $\beta_1^i$. Consider an expansion $\sum_{m\ge0}c_mF(t)^m$ with
$|c_m|\le1$. For $m\ge1$, the coefficient of $t^j$ in
$F(t)^m$ is
\[
\sum_{\substack{i_1+\cdots+i_m=j\\i_1,\ldots,i_m\ge1}}
F_{i_1}\cdots F_{i_m}.
\]
It is easy to see that $m\le j$. For each ordered
decomposition, the total weight of the product is at most
$\beta_1^{i_1}\cdots \beta_1^{i_m}=\beta_1^j$. There are $\binom{j-1}{m-1}$
such decompositions. Summing over $m$ gives the bound
\[
\beta_1^j\sum_{m=1}^j\binom{j-1}{m-1}
=2^{j-1}\beta_1^j\le(2\beta_1)^j.
\]
This rule applies to the geometric series, for which $c_m=1$,
and to the inverse square root, for which
$c_m=\binom{-1/2}{m}
=\frac{(-1/2)(-1/2-1)\cdots(-1/2-m+1)}{m!}$ and $|c_m|\le1$.

\medskip
\noindent\textbf{Application to the coefficients of $L(tV)$.} Let $b_j$ denote the total weight of the terms in $L_j$. We apply
the two rules in the order in which the factors of
$L(tV)=W(tV)^*T_{\Omega}(tV)W(tV)$ are constructed.

\textup{(i)} \noindent\emph{The coefficients of $T_{\Omega}(tV)$.}
First apply Rule~2 to $F(t)=A(tV)-\operatorname{Id}_{\mathcal E}
=\sum_{i\ge1}t^iA_i$. Since
\[
A(tV)^{-1/2}
=(\operatorname{Id}_{\mathcal E}+F(t))^{-1/2}
=\sum_{m\ge0}\binom{-1/2}{m}F(t)^m,
\]
and the coefficient $A_i$ has weight at most $C^i$,
Rule~2 gives a total weight of at most $(2C)^i$
for the coefficient of degree $i\ge1$ in $A(tV)^{-1/2}$. Next apply the three-factor case of
Rule~1 to
\[
T_{\Omega}(tV)
=A(tV)^{-1/2}B(tV)A(tV)^{-1/2}.
\]
It follows that the total weight of $T_{\Omega,i}$ is at most $K_1^i$
for some fixed $K_1\ge1$.

\textup{(ii)} \noindent\emph{The coefficients of $P(tV)$.}
Set $D(t)=T_{\Omega}(tV)-B_0$ as in Step~3. For each
$z\in\Gamma$, we apply Rule~2 to $F(t)=R_0(z)D(t)$ in the
resolvent expansion
\[
(z\operatorname{Id}_{\mathcal E}-T_{\Omega}(tV))^{-1}
=\sum_{m\ge0}\bigl(R_0(z)D(t)\bigr)^mR_0(z).
\]
Indeed, if we let $\beta_0=\max\{1,\sup_{z\in\Gamma}\|R_0(z)\|\}$,
the total weight of the coefficient $R_0(z)T_{\Omega,i}$ is at most
$\beta_0 K_1^i\le(\beta_0 K_1)^i$. Thus Rule~2 gives a bound uniform in $z$.
The final factor $R_0(z)$ and the contour integral
multiply this bound by at most $\beta_0|\Gamma|/(2\pi)$. Absorbing
this constant into the base of the exponential gives a
bound $K_2^i$ for the total weight of $P_i$.

\textup{(iii)} \noindent\emph{The coefficients of $W(tV)$.}
Since $N_i=P_0P_iP_0|_{E_0}$ and $P_0$ is an orthogonal
projection, the same bound holds for $N_i$. Thus, we apply Rule~2
to $F(t)=N(tV)-\operatorname{Id}_{E_0}$ to obtain an
exponential bound for the coefficients of $N(tV)^{-1/2}$.
Then we use the two-factor case of Rule~1 to $W(tV)=\bigl(P(tV)|_{E_0}\bigr)N(tV)^{-1/2}$. This gives an exponential bound for the total weight of $W_i$.
Taking adjoints gives the same bound for $W_i^*$.

\textup{(iv)} \noindent\emph{The coefficients of $L(tV)$.}
Finally, it is easy to see that the coefficient of degree $j$ in $L(tV)=W(tV)^*T_{\Omega}(tV)W(tV)$
is
\[
L_j=\sum_{\substack{a+b+c=j\\a,b,c\ge0}}
W_a^*T_{\Omega,b}W_c.
\]
The three-factor case of Rule~1 now
gives $b_j\le K_3^j, j\ge1$, for constant $K_3$ independent of $V$.

\medskip
\noindent\textbf{From $L(tV)$ to its inverse.}
The expansion of $L(tV)^{-1}$ gives, for $k\ge1$,
\[
H_{\Omega,k}(V,\ldots,V)
=\sum_{m=1}^k(-1)^m\lambda_0^{m+1}
\sum_{\substack{i_1+\cdots+i_m=k\\i_1,\ldots,i_m\ge1}}
L_{i_1}\cdots L_{i_m}.
\]
For each ordered decomposition, the total weight of the expanded
product is at most
$b_{i_1}\cdots b_{i_m}\le K_3^{i_1+\cdots+i_m}=K_3^k$.
Counting the decompositions of length $m$ as in Rule~2 and
retaining the scalar factor $\lambda_0^{m+1}$, we obtain
\[
a_k
\le K_3^k\sum_{m=1}^k\binom{k-1}{m-1}\lambda_0^{m+1}
=K_3^k\lambda_0^2(1+\lambda_0)^{k-1},
\qquad k\ge2.
\]
Since $\lambda_0>0$ is fixed, there exists a constant $D\ge1$,
independent of $k$ and $V$, such that $a_k\le D^k$ for all $k\ge2$.

\bigskip
\noindent\textbf{Conflict of interest and data availability statements.}
\par\medskip
\indent The authors declare no conflict of interest.
This article presents theoretical mathematical results.
No datasets were generated or analyzed, so data sharing is not applicable.

\bigskip
\noindent\textbf{Declaration of AI use.}
\par\medskip
\indent During the preparation of this work, the authors used ChatGPT
and Doubao to improve the language and grammar and to assist with
checking some calculations.
The authors carefully reviewed and edited the content as needed.
They take full responsibility for the final version of this manuscript.


\begin{thebibliography}{99}

\bibitem{AF}
P. R. S. Antunes and P. Freitas,
Numerical optimization of low eigenvalues of the Dirichlet and Neumann Laplacians,
\emph{J. Optim. Theory Appl.} \textbf{154} (2012), 235--257.

\bibitem{AO}
P. R. S. Antunes and E. Oudet,
Numerical results for extremal problem for eigenvalues of the Laplacian,
in A. Henrot (ed.), \emph{Shape Optimization and Spectral Theory},
De Gruyter Open, Warsaw, 2017, Chapter~11, 398--411.

\bibitem{AB}
M. S. Ashbaugh and R. D. Benguria,
A sharp bound for the ratio of the first two eigenvalues of Dirichlet Laplacians and extensions,
\emph{Ann. of Math.} (2) \textbf{135} (1992), 601--628.

\bibitem{ABR}
S. Axler, P. Bourdon, and W. Ramey,
\emph{Harmonic Function Theory}, 2nd ed., Graduate Texts in Mathematics~137,
Springer, New York, 2001.

\bibitem{Bossel}
M.-H. Bossel,
Membranes \'elastiquement li\'ees: extension du th\'eor\`eme de Rayleigh--Faber--Krahn et de l'in\'egalit\'e de Cheeger,
\emph{C. R. Acad. Sci. Paris S\'er. I Math.} \textbf{302} (1986), 47--50.

\bibitem{BDPV}
L. Brasco, G. De Philippis, and B. Velichkov,
Faber--Krahn inequalities in sharp quantitative form,
\emph{Duke Math. J.} \textbf{164} (2015), 1777--1831.
\href{https://arxiv.org/abs/1306.0392}{arXiv:1306.0392}.

\bibitem{Bucur}
D. Bucur,
Minimization of the $k$-th eigenvalue of the Dirichlet Laplacian,
\emph{Arch. Ration. Mech. Anal.} \textbf{206} (2012), 1073--1083.

\bibitem{DL}
M. Dambrine and J. Lamboley,
Stability in shape optimization with second variation,
\emph{J. Differential Equations} \textbf{267} (2019), 3009--3045.
Corrected author version, 2021,
\href{https://webusers.imj-prg.fr/~jimmy.lamboley/fichiers/StabilityCorrectedFinal.pdf}{available online}.

\bibitem{DanersRobin}
D. Daners,
A Faber--Krahn inequality for Robin problems in any space dimension,
\emph{Math. Ann.} \textbf{335} (2006), 767--785.

\bibitem{Davies}
E. B. Davies,
\emph{Spectral Theory and Differential Operators},
Cambridge Studies in Advanced Mathematics~42,
Cambridge University Press, Cambridge, 1995.

\bibitem{ESI}
A. El Soufi and S. Ilias,
Domain deformations and eigenvalues of the Dirichlet Laplacian in a Riemannian manifold,
\emph{Illinois J. Math.} \textbf{51} (2007), 645--666.

\bibitem{Faber}
G. Faber,
Beweis, da\ss{} unter allen homogenen Membranen von gleicher Fl\"ache und gleicher Spannung die kreisf\"ormige den tiefsten Grundton gibt,
\emph{Sitzungsber. Bayer. Akad. Wiss. Math.-Phys. Kl.} (1923), 169--172.
\href{https://publikationen.badw.de/de/003399311}{Original article}.

\bibitem{FKRobin}
P. Freitas and D. Krej\v{c}i\v{r}\'{i}k,
The first Robin eigenvalue with negative boundary parameter,
\emph{Adv. Math.} \textbf{280} (2015), 322--339.
\href{https://arxiv.org/abs/1403.6666}{arXiv:1403.6666}.

\bibitem{GT}
D. Gilbarg and N. S. Trudinger,
\emph{Elliptic Partial Differential Equations of Second Order},
2nd ed., Classics in Mathematics, Springer, Berlin, 2001.

\bibitem{Henrot}
A. Henrot,
\emph{Extremum Problems for Eigenvalues of Elliptic Operators},
Frontiers in Mathematics, Birkh\"auser, Basel, 2006.

\bibitem{HP}
A. Henrot and M. Pierre,
\emph{Shape Variation and Optimization: A Geometrical Analysis},
EMS Tracts in Mathematics~28, European Mathematical Society, Z\"urich, 2018.
\href{https://ems.press/books/etm/149}{Publisher's page}.

\bibitem{HS}
L. Hermi and N. Saito,
On Rayleigh-type formulas for a non-local boundary value problem
associated with an integral operator commuting with the Laplacian,
\emph{Appl. Comput. Harmon. Anal.} \textbf{45} (2018), 59--83.

\bibitem{Kato}
T. Kato,
\emph{Perturbation Theory for Linear Operators},
reprint of the 1980 edition, Classics in Mathematics,
Springer, Berlin, 1995.

\bibitem{KennedyRobin}
J. B. Kennedy,
An isoperimetric inequality for the second eigenvalue of the Laplacian with Robin boundary conditions,
\emph{Proc. Amer. Math. Soc.} \textbf{137} (2009), 627--633.

\bibitem{KennedyHigher}
J. B. Kennedy,
On the isoperimetric problem for the higher eigenvalues of the Robin and Wentzell Laplacians,
\emph{Z. Angew. Math. Phys.} \textbf{61} (2010), 781--792.

\bibitem{Krahn25}
E. Krahn,
\"Uber eine von Rayleigh formulierte Minimaleigenschaft des Kreises,
\emph{Math. Ann.} \textbf{94} (1925), 97--100.

\bibitem{Krahn26}
E. Krahn,
\"Uber Minimaleigenschaften der Kugel in drei und mehr Dimensionen,
\emph{Acta Comment. Univ. Dorpat. A} \textbf{9} (1926), no.~1, 1--44.

\bibitem{Lorch}
L. Lorch,
Some inequalities for the first positive zeros of Bessel functions,
\emph{SIAM J. Math. Anal.} \textbf{24} (1993), 814--823.

\bibitem{McLean}
W. McLean,
\emph{Strongly Elliptic Systems and Boundary Integral Equations},
Cambridge University Press, Cambridge, 2000.

\bibitem{DLMF}
National Institute of Standards and Technology,
\emph{NIST Digital Library of Mathematical Functions},
Chapter~10, Bessel functions, especially Sections~10.2, 10.6, 10.19, and 10.21.
\url{https://dlmf.nist.gov/10}.

\bibitem{OK}
B. Osting and C.-Y. Kao,
Minimal convex combinations of three sequential Laplace--Dirichlet eigenvalues,
\emph{Appl. Math. Optim.} \textbf{69} (2014), 123--139.

\bibitem{Oudet}
E. Oudet,
Numerical minimization of eigenmodes of a membrane with respect to the domain,
\emph{ESAIM Control Optim. Calc. Var.} \textbf{10} (2004), 315--330.

\bibitem{PS}
G. P\'olya and G. Szeg\H{o},
\emph{Isoperimetric Inequalities in Mathematical Physics},
Annals of Mathematics Studies~27,
Princeton University Press, Princeton, NJ, 1951.

\bibitem{PTZ}
Y. Privat, E. Tr\'elat, and E. Zuazua,
Spectral shape optimization for the Neumann traces of the
Dirichlet-Laplacian eigenfunctions,
Calc. Var. Partial Differential Equations \textbf{58} (2019),
Article 64.

\bibitem{Rayleigh}
Lord Rayleigh,
\emph{The Theory of Sound}, Vol.~I, 2nd ed.,
Macmillan, London, 1894.

\bibitem{SZ}
J. Soko\l{}owski and J.-P. Zol\'esio,
\emph{Introduction to Shape Optimization: Shape Sensitivity Analysis},
Springer Series in Computational Mathematics~16,
Springer, Berlin, 1992.

\bibitem{Szego}
G. Szeg\H{o},
Inequalities for certain eigenvalues of a membrane of given area,
\emph{J. Rational Mech. Anal.} \textbf{3} (1954), 343--356.

\bibitem{Watson}
G. N. Watson,
\emph{A Treatise on the Theory of Bessel Functions},
2nd ed., Cambridge University Press, Cambridge, 1944.

\bibitem{Weinberger}
H. F. Weinberger,
An isoperimetric inequality for the $N$-dimensional free membrane problem,
\emph{J. Rational Mech. Anal.} \textbf{5} (1956), 633--636.

\bibitem{WK}
S. A. Wolf and J. B. Keller,
Range of the first two eigenvalues of the Laplacian,
\emph{Proc. Roy. Soc. London Ser. A} \textbf{447} (1994), 397--412.

\end{thebibliography}
\end{document}